\documentclass[11pt]{article}
\usepackage{amsmath} 
\usepackage{amsthm} 
\usepackage{amssymb}
\usepackage{enumerate}

\def\today{\number\day\space
\ifcase\month\or
  January\or February\or March\or April\or May\or June\or
  July\or August\or September\or October\or November\or 
  December%
\fi
\space\number\year}

\swapnumbers
\newtheorem{brslt}{}[section]
\newcommand{\start}[2]{\begin{brslt} \label{#2}  #1}
\newcommand{\trats}{\end{brslt}}

\newcommand{\startnm}[1]{\begin{enumerate}[\upshape(#1)]}
\newcommand{\tratsnm}{\end{enumerate}}

\newcommand{\alrt}[1]{{\em #1}}

\def\multdef#1#2{\actdef{#1}{#2}}
\def\actdef#1#2#3{%
   \if #3. \def\next{}%
   \else\expandafter\def\csname #1#3\endcsname{{#2{#3}}}%
        \def\next{\multdef{#1}{#2}}
   \fi%
 \next}

\multdef{cal}{\mathcal}
{A}{B}{C}{D}{E}{F}{G}{H}{I}{J}{K}{L}{M}
{N}{O}{P}{Q}{R}{S}{T}{U}{V}{W}{X}{Y}{Z}.

\multdef{bb}{\mathbb}
{A}{B}{C}{D}{E}{F}{G}{H}{I}{J}{K}{L}{M}
{N}{O}{P}{Q}{R}{S}{T}{U}{V}{W}{X}{Y}{Z}.

\newcommand{\goesto}{\longrightarrow}
\newcommand{\iso}{\simeq}

\renewcommand{\phi}{\varphi}
\newcommand{\eps}{\epsilon}

\newcommand{\pf}{{\scshape Proof.} }
\newcommand{\thm}{{\scshape Theorem.} }
\newcommand{\lem}{{\scshape Lemma.} }
\newcommand{\prop}{{\scshape Proposition.} }
\newcommand{\cor}{{\scshape Corollary.} }

\newcommand{\rmks}{{\scshape Remarks.} }
\newcommand{\hyp}{{\scshape Hypothesis.} }

\multdef{name}{\operatorname}
    {Sym}{Aut}{SL}.
\newcommand{\Sym}[1]{\nameSym{(#1)}}
\newcommand{\Aut}[2][]{\nameAut_{#1}{(#2)}}
\newcommand{\SL}[2]{\nameSL_{#1}{(#2)}}

\newcommand{\set}[2]{\{\, #1 \mid #2\, \}}
\newcommand{\dspace}[1]{\vspace*{#1\baselineskip}}
\newcommand{\done}{{\qquad\hfill\ensuremath{\Box}}}
\newcommand{\zt}{{\bbZ/2\bbZ}}
\newcommand{\lam}{\lambda}
\newcommand{\Lam}{\Lambda}

\newcommand{\oA}{\operatorname{1A}}
\newcommand{\tA}{\operatorname{2A}}
\newcommand{\tB}{\operatorname{2B}}
\newcommand{\tC}{\operatorname{3C}}

\newcommand{\ad}{\operatorname{ad}}
\newcommand{\VOA}{VOA}

\newcommand{\sym}{J}
\newcommand{\Vp}[2]{{#1}^{\sym}(#2)}
\newcommand{\Clif}[1]{\operatorname{Cl}(#1)}
\newcommand{\SClif}[2][\sym]{\operatorname{Cl}^{#1}(#2)}
\newcommand{\bb}[1]{b_{#1}}
\newcommand{\onehalf}{\nfrac{1}{2}}
\newcommand{\halfeta}{\nfrac{\eta}{2}}
\newcommand{\axisset}{\calA}
\newcommand{\Matsuo}[1]{\operatorname{M}(#1)}
\newcommand{\fusion}[1][F]{\mathfrak{#1}}

\newcommand\lla{{\langle\!\langle}}
\newcommand\rra{{\rangle\!\rangle}}

\newcommand{\nfrac}[2]{\textstyle{\frac{#1}{#2}}}
\newcommand{\name}[1]{{\sf(#1)}}

\newcommand{\nm}[1]{{\rm(#1)}}

\title{Primitive axial algebras of Jordan type}

\author{
J.I. Hall\\
Department of Mathematics\\
Michigan State University\\
Wells Hall\\
619 Red Cedar Road\\
East Lansing, Michigan 48840\\
United States of America\\
{\tt jhall@math.msu.edu}
\\[5pt]
F. Rehren\\
School of Mathematics\\
Watson Building\\
University of Birmingham\\
Edgbaston\\
Birmingham, B15 2TT\\
United Kingdom\\
{\tt rehrenf@maths.bham.ac.uk}
\\[5pt]
S. Shpectorov\\
School of Mathematics\\
Watson Building\\
University of Birmingham\\
Edgbaston\\
Birmingham, B15 2TT0\\
United Kingdom\\
{\tt s.shpectorov@bham.ac.uk}
}

\date{Version of: 7 March 2014}

\begin{document}
\maketitle

\begin{abstract}

An \alrt{axial algebra} over the field $\bbF$ is a commutative algebra
generated by idempotents whose adjoint action has multiplicity-free
minimal polynomial. For semisimple associative algebras this leads to
sums of copies of $\bbF$.  Here we consider the first nonassociative
case, where adjoint minimal polynomials divide $(x-1)x(x-\eta)$ for
fixed $0\neq\eta\neq 1$. Jordan algebras arise when $\eta=\onehalf$,
but our motivating examples are certain Griess algebras of vertex
operator algebras and the related Majorana algebras. We study a class
of algebras, including these, for which axial automorphisms like those
defined by Miyamoto exist, and there classify the $2$-generated
examples.  For $\eta \neq \onehalf$ this implies that the Miyamoto
involutions are $3$-transpositions, leading to a classification.
\end{abstract}

\section{Introduction}
\label{sec-intro}

Throughout we consider commutative $\bbF$-algebras $A$ where $\bbF$ is
a field of characteristic not equal to two.  We emphasize that our
algebras will usually be nonassociative and may not have an identity
element.

For the element $a$ of $A$ and $\lam \in \bbF$, the $\lam$-eigenspace
for the adjoint $\bbF$-endomorphism $\ad_a$ of $A$ will be denoted
$A_\lam(a)$ (where we allow the possibility $A_\lam(a)=0$).  If $A$ is
an associative algebra and $a$ is an idempotent element, then
$A=A_1(a)\oplus A_0(a)$---the adjoint of the idempotent is semisimple
with at most the two eigenvalues $0$ and $1$.  Here we are interested
in the minimal nonassociative case---semisimple idempotents whose
adjoint eigenvalues are drawn from the set $\Lam = \{1,0,\eta\}$ for
some $\eta \in \bbF$ with $0 \neq \eta \neq 1$.

An idempotent whose adjoint is semisimple will be called an
\alrt{axis}. A commutative algebra generated by axes is then an
\alrt{axial algebra}.  The commutative algebra $A$ over $\bbF$ (not of
characteristic two) is a \alrt{primitive axial algebra of Jordan type
$\eta$} provided it is generated by a set of axes with each member
$a$ satisfying: 
\startnm{a}
\item  
$A=A_1(a) \oplus A_0(a) \oplus A_\eta(a)$.
\item 
$A_1(a)=\bbF a$. 
\item
$A_0(a)$ is a subalgebra of $A$.
\item
For all $\delta,\eps \in \pm$,
\[
A_\delta(a) A_\eps(a) \subseteq A_{\delta\eps}(a)\,,
\]
where $A_+(a) = A_1(a)\oplus  A_{0}(a)$
and $A_-(a) = A_{\eta}(a)$.
\tratsnm

Examples include Jordan algebras that are generated by idempotents
\cite{jac-jordan}. These occur for $\eta=\frac{1}{2}$, although this
is the case in which we say the least. Instead our motivation comes
from the values $\eta=\frac{1}{4}$ and $\eta=\frac{1}{32}$, which
arise as special cases of $\Lam = \{1,0,\frac{1}{4},\frac{1}{32}\}$.
Algebras of this latter type are provided by Griess algebras
associated with vertex operator algebras and Majorana algebras
\cite{iv-book,Matsuo,Miyamoto,S}.

A major accomplishment in the Griess algebra case was Sakuma's Theorem
\cite{S} which classified all $2$-generated subalgebras. See also
\cite{iv-book,ipss,RS}.  The following similar theorem is a central
result of this paper.

\start{\thm}{thm-2gen-jordan-field}
Let $\bbF$ be a field of characteristic not two with $\eta \in \bbF$
for $0 \neq \eta \neq 1$.  Let $A$ be a primitive axial $\bbF$-algebra
of Jordan type $\eta$ that is generated by two axes.  Then we have one
of the following:
\startnm{1}
\item $A$ is an algebra $\bbF$ of type $\oA$ over $\bbF$;
\item $A$ is an algebra $\bbF \oplus \bbF$ of type $\tB$ over $\bbF$;
\item $A$ is an algebra of type $\tC(\eta)$ of dimension $3$ over
  $\bbF$;
\item $\eta=-1$ and $A$ is an algebra of type $\tC(-1)^*$ of dimension
  $2$ over $\bbF$;
\item $\eta=\nfrac{1}{2}$ and $A$ is isomorphic to the $3$-dimensional
  symmetric Jordan Clifford algebra $\SClif{\bbF^2,\bb{\delta}}$,
  where the symmetric bilinear form $\bb{\delta}$ on $\bbF^2$ is given
  by $\bb{\delta}(v_i,v_i)=2$ and $\bb{\delta}(v_0,v_1)=\delta\neq 2$
  for its basis $v_0,v_1$.
\item $\eta=\nfrac{1}{2}$ and $A$ is isomorphic to the $2$-dimensional
  special Jordan algebra $\SClif[0]{\bbF^2,\bb{2}}$ or the
  $3$-dimensional Jordan algebra $\SClif[00]{\bbF^2,\bb{2}}$, where
  the degenerate symmetric bilinear form $\bb{2}$ on $\bbF^2$ is given
  by $\bb{2}(v_i,v_j)=2$ for its basis $v_0,v_1$.
\tratsnm
\trats

\noindent
For the definitions and discussion of the various examples, see
Section \ref{sec-eg}.

The restriction \nm{d} provides a $\zt$-grading of $A$ for each axis
$a$. Equivalently, the linear transformation of $A$ that acts as the
identity on $A_+(a)$ and negates everything in $A_-(a)$ is an
automorphism of $A$.  The resulting automorphisms of order $2$ will be
called \alrt{Miyamoto involutions}, since in the Griess algebra
context they were first noticed and used to great effect by Miyamoto
in \cite{Miyamoto}.  Sakuma's proof and our proof of Theorem
\ref{thm-2gen-jordan-field} make critical use of the dihedral group
generated by the Miyamoto involutions corresponding to the two
generators.

If $a$ is an axis and $g$
is an automorphism of $A$, then $a^g$ is also
an axis.
For a generating set $\axisset$ of axes, let 
$\bar{\axisset}$ be the smallest set of
axes with the properties:
\startnm{i}
\item $\axisset \subseteq \bar{\axisset}$.
\item If $b \in \bar{\axisset}$ and $\tau$ is the Miyamoto involution
  associated with $b$, then $\bar{\axisset}^\tau \subseteq
  \bar{\axisset}$. 
\tratsnm

As a consequence of the theorem, every product of two members of
$\axisset$, and indeed of $\bar{\axisset}$, is in the $\bbF$-span of
$\bar{\axisset}$. Therefore

\start{\cor}{cor-2gen}
Let $A$ be an axial algebra of Jordan type $\eta$ over a field $\bbF$
of characteristic not two that is generated by the set $\axisset$ of
axes.  Then $A$ is spanned as $\bbF$-space by the axes of
$\bar{\axisset}$.
\trats

The theorem readily leads to

\start{\thm}{thm-miyamoto-neq2}
Let $A$ be an axial algebra of Jordan type $\eta \neq \frac{1}{2}$
over a field of characteristic not two that is generated by the set
$\axisset$ of axes.  Then the Miyamoto involutions corresponding to
$\bar{\axisset}$ form a normal set of $3$-transpositions in the
automorphism group of $A$ that they generate.
\trats

\noindent
For the definition and discussion of $3$-transpositions, see Section
\ref{sec-auto}.

The theorem and results from \cite{CuHa} imply that in finitely
generated algebras with $\eta\neq \frac{1}{2}$, the set
$\bar{\axisset}$ is finite. Together with Corollary \ref{cor-2gen}
this leads to

\start{\cor}{cor-neq2-loc-fin}
If the axial algebra $A$ of Jordan type $\eta \neq \frac{1}{2}$ over
the field $\bbF$ of characteristic not two is finitely generated then
it is finite dimensional as a vector space over $\bbF$.
\trats

The theorem is in general false for Jordan type $\eta=\frac{1}{2}$,
and we suspect that the corollary is also false in that case.
Certainly in that case a finitely generated algebra can have
$\bar{\axisset}$ infinite; see remarks near the beginning of Section
\ref{sec-auto}.

We have a converse to Theorem \ref{thm-miyamoto-neq2}.

\start{\thm}{thm-matsuo-bis}
Let $D$ be a normal set of $3$-transpositions in the group $G=\langle
D \rangle$. For any field $\bbF$ not of characteristic two and any
$\eta \in \bbF$ with $0\neq\eta\neq 1$, the space $M=\bbF D$ can be
given the structure of a primitive axial algebra of Jordan type $\eta$
on which the elements of $D$ act as Miyamoto involutions.  The algebra
$M$ admits a nonzero symmetric and associative bilinear form
\(
\lla \cdot,\cdot \rra \colon A \times A \goesto \bbF\,.
\)
\trats

This is proven in a more precise form in Theorem \ref{thm-matsuo} and
Corollary \ref{cor-matsuo-frob} below.

\dspace{.5}

The historical context for the commutative algebras discussed in this
paper has three separate branches.  The first, discussed above, views
axial algebras of Jordan type $\eta$ as a first step away from
semisimple, associative algebras---a step far enough away to include
new and interesting examples, such as the Jordan algebras generated by
idempotents, but not so far as to defy meaningful classification.

The second branch, also mentioned above, was the actual motivation for
the present work.  In the early 1970's Bernd Fischer and Robert Griess
independently found evidence for the existence of the Monster sporadic
simple group $\bbM$. Soon after, it was noted that the smallest
faithful $\bbR$-module for $\bbM$ might well have dimension $196883$,
and Simon Norton observed that such a module would admit a
commutative, nonassociative algebra structure \cite{Gr76,CoNo79}.

Bob Griess \cite{Gr81} constructed this algebra and hence $\bbM$ as an
automorphism group (by hand). In the full treatment \cite{griess} he
preferred an algebra of dimension $196884$, including a trivial
$\bbM$-submodule. Conway \cite{conway} used a deformation $B^\natural$
of the $196884$ algebra to give a new construction. He noted the
existence in $B^\natural$ of idempotents (after appropriate scaling),
one for each $\tA$ involution of the Monster, with adjoint minimal
polynomial $(x-1)x(x-\nfrac{1}{4}) (x-\nfrac{1}{32})$; he called these
idempotents \alrt{axial vectors}.

Motivated in part by the ``Monstrous Moonshine'' conjectures
\cite{CoNo79}, Borcherds \cite{Borcherds} codified vertex operators,
and Frenkel, Lepowsky, and Meurman \cite{FLM88} constructed a vertex
operator algebra $V^\natural$ whose graded piece $V^\natural_2$
inherits from the \VOA\ a natural commutative algebra structure
isomorphic to $B^\natural$.  The algebra $V^\natural$ belongs to a
large class of \VOA\ for which $V_2$ always admits a natural structure
as commutative algebra.  These commutative algebras are the
\alrt{Griess algebras}, and Miyamoto \cite{Miyamoto} observed that in
them each conformal vector of central charge $\onehalf$ can be viewed
as an axis, in the sense that there is an involutory automorphism
acting in a prescribed way relative to the
$\{1,0,\nfrac{1}{4},\nfrac{1}{32}\}$-eigenspaces of each of these
conformal vectors. This effectively reverses Conway's construction of
axes from the $2A$ involutions of $V^\natural_2=B^\natural$.
Miyamoto, Kitazume, and others \cite{Miyamoto,KiMi01} then studied the
possible groups that can be generated by these \alrt{Miyamoto
  involutions} within the automorphism groups of Griess algebras and
their associated \VOA s. Of particular importance for this paper is
the work of Sakuma \cite{S}, which described all groups generated by
two such involutions, and that of Matsuo \cite{Matsuo}, which
completed the study of the case where only the eigenvalues
$\{1,0,\nfrac{1}{4}\}$ occur. Indeed the original version
\cite{Matsuo-ver1} of \cite{Matsuo} discussed the more general case of
Griess algebras with axes admitting only three eigenvalues and noted
there that the Miyamoto involutions are $3$-transpositions, an
observation due to Miyamoto \cite[Theorem~6.13]{Miyamoto} in the
$\{1,0,\nfrac{1}{4}\}$ case.  Matsuo's unpublished original thus
contains versions of several of the main results of this paper, albeit
in the more restricted context of Griess algebras.

In an effort to divorce the properties of Griess algebras from the
\VOA s that envelope them, Ivanov \cite{iv-book} introduced
\alrt{Majorana algebras}. From our point of view (see Section
\ref{sec-fus}) these, and so especially the Griess algebra examples,
are real, Frobenius, primitive axial algebras with fusion table
\renewcommand{\arraystretch}{1.75}
\begin{center}
\begin{tabular}{c||ccc|c}
$\star$ & 1  & 0 & $\nfrac{1}{4}$ & $\nfrac{1}{32}$\\
\hline\hline
1 & 1  & $\emptyset$ & $\nfrac{1}{4}$ & $\nfrac{1}{32}$\\
0 &$\emptyset$& 0 & $\nfrac{1}{4}$ & $\nfrac{1}{32}$\\
$\nfrac{1}{4}$ &
$\nfrac{1}{4}$ & $\nfrac{1}{4}$ & 1,0 & $\nfrac{1}{32}$\\
\hline
$\nfrac{1}{32}$ &
$\nfrac{1}{32}$ & $\nfrac{1}{32}$ & 
$\nfrac{1}{32}$ & 1,0,$\nfrac{1}{4}$\\
\end{tabular}
\end{center}
Ivanov, Pasechnik, Seress, and Shpectorov \cite{ipss} extended
Sakuma's $2$-generator theorem to Majorana algebras.

The third and final contextual branch for this paper begins with Simon
Norton \cite{No75,No88}, who constructed a commutative, nonassociative
algebra with automorphism group a triple cover of the sporadic group
$Fi_{24}$. Although the term ``Norton algebra'' remains loosely
defined \cite{Sm77}, such constructions start with a partial linear
space of order $2$ (see Section \ref{sec-matsuo}) as the skeleton for
a presentation of a commutative algebra via its $2$-generated
subalgebras.  Work along these lines was done by many, in part to
characterize various groups as the automorphism groups of algebra
structures on related ``natural'' modules.  A case in point is the
paper \cite{Har84} in which, among other things, Harada constructed a
$6$-dimensional commutative nonassociative algebra with automorphism
group $\SL{3}{2}$, using as skeleton $PG(2,2)$---the projective plane
of order $2$.  The unpublished \cite{MaMa} of Matsuo and Matsuo
studied in detail a related $1$-parameter family of complex
$7$-dimensional algebras, each with skeleton $PG(2,2)$, generically
having $\SL{3}{2}$ as automorphism group. The unpublished original
\cite{Matsuo-ver1} pursued and generalized this construction, noting
that for Griess algebras the appropriate partial linear spaces to
consider are the Fischer spaces---the geometric counterparts to
$3$-transposition groups.  The algebras of \cite{MaMa} and
\cite{Matsuo-ver1} are examples of the \alrt{Matsuo algebras}
introduced here in Section \ref{sec-matsuo}.  Those corresponding to
Fischer spaces were also discussed in \cite{Re13}.

One possible path for generalization is to replace the idempotence
condition by the requirement that $a^2= k a$, for some constant
$k$. If $k$ is a nonzero field element then rescaling gets us back to
the idempotent case, but $a^2=0$ can lead to interesting algebras. Let
$\Matsuo{\Pi,\onehalf,\bbQ}$ be the rational Matsuo algebra for the
partial linear space $\Pi$ of order $2$, and let $A$ be its
$\bbZ$-subalgebra spanned the various elements $2a_p$, for $p$ a point
of $\Pi$. Then $A/2A$ is a commutative, nonassociative
$\bbF_2$-algebra generated by elements with square $0$. For various
choices of $\Pi$ these algebras and their quotients give rise to Lie
algebras and ``near-Lie'' algebras that have interesting automorphism
groups, as can be seen in \cite{Cu05} and \cite{Jos12}.

\dspace{.5}

We conclude with a brief outline of the paper. After this introduction
we present the basic definitions for axial algebras and their related
fusion rules. The fundamental objects are the semisimple idempotents,
and we include a characterization of semisimple associative
commutative algebras in this context. Of particular interest are the
algebras for which the fusion rules guarantee the existence of
Miyamoto involutory automorphisms.

In the algebras of Jordan type the idempotents have at most three
adjoint eigenvalues---$1$ and $0$ (as expected) and some $\eta \ne
0,1$.  Section \ref{sec-eg} describes various examples of primitive
axial algebras of Jordan type. The focus is on $2$-generated algebras,
but the special case of Jordan algebras is discussed in greater
generality.

Section \ref{sec-poc} contains the proof of Theorem
\ref{thm-2gen-jordan-field} of Sakuma type. The next section then
describes the impact of the Sakuma theorem on the automorphism group
of an arbitrary primitive axial algebra of Jordan type
$\eta\neq\onehalf$. Specifically the Miyamoto involutions provide a
normal set of $3$-transpositions, as detailed in Theorem
\ref{thm-miyamoto-neq2}.  The classifications of \cite{fischer} and
\cite{CuHa} are then available.

In Section \ref{sec-matsuo} we provide a constructional converse to
the results of the previous section.  Beginning with a partial linear
space of order $2$, an algebra is constructed---the \alrt{Matsuo
  algebra}.  In particular it is proven, as in Theorem
\ref{thm-matsuo-bis}, that every normal set of $3$-transpositions
arises as a set of Miyamoto involutions for an appropriate primitive
axial algebra of Jordan type $\eta$, for arbitrary $\eta$. Indeed the
Matsuo algebras that have Jordan type $\eta$ are precisely those with
Fischer spaces as skeletons; see Theorem \ref{thm-felix}.

Griess algebras come equipped with an associative form, and in the
final Section \ref{sec-frob} we describe the circumstances under which
a Matsuo algebra admits such a form. This includes all the axial
algebras of Theorem \ref{thm-matsuo-bis}.

\dspace{.5}

Parts of the present article have much in common with the unpublished
paper of Matsuo \cite{Matsuo-ver1} (although he makes some assumptions
not made here---he assumes characteristic $0$ and proves a version of
Theorem \ref{thm-2gen-jordan-field} only for Griess algebras). We
warmly thank Professor Matsuo for pointing out \cite{Matsuo-ver1} and
providing us with hard copy of the otherwise unavailable~\cite{MaMa}.

\section{Fusion and axial algebras}
\label{sec-fus}
\subsection{Decomposition, fusion, and grading}

For the $\bbF$-algebra $A$ and set $I$, a \alrt{decomposition} results
from writing $A$ as a direct sum of subspaces indexed by $I$:
\[
A = \bigoplus_{i \in I} A_i\,.
\]
Each of our decompositions will be accompanied by a \alrt{fusion rule}
\[
\fusion\colon I \times I \goesto 2^I\,,
\]
a map that takes each ordered pair $i,j$ of indices to a member $i
\star j$ of the power set of $I$ and that encodes the \alrt{fusion}
information
\[
A_iA_j \subseteq 
 \bigoplus_{k \in i\star j} A_k\,. 
\]
We then have an \alrt{$\fusion$-decomposition}.

Of course every decomposition admits the trivial fusion rule with $i
\star j=I$, for all $i,j \in I$.  Equally well every decomposition has
a unique \alrt{minimal} fusion rule, where each $i \star j$ is chosen
with cardinality as small as possible.

We are interested in rules where all $i\star j$ have small
cardinality, in which case the fusion rule may be easily presented in
a \alrt{fusion table}.  In particular, if each $i \star j$ has
cardinality $1$, then the fusion table gives the multiplication table
for a magma $(I,\star)$, and the encoded fusion properties describe a
\alrt{grading} of $A$ by that magma.

The algebras we study will be graded by $\zt$, the cyclic group of
order two.  We there take $I =\{\pm 1\}=\{\pm\}$ and $i\star j = ij$.
It is well-known and easy to see that in characteristic not two the
existence of a $\zt$-grading is equivalent to the existence of an
automorphism of order (at most) two:

\start{\prop}{prop-automorphism}
Let $A$ be an $\bbF$-algebra with the characteristic of $\bbF$ not
two.  Suppose $A= A_+ \oplus A_-$. Then the following are equivalent:
\startnm{1}
\item 
The indices give a $\zt$-grading of $A$.
\item
$A$ has an automorphism of order at most two that fixes each element
  of $A_+$ and takes each element of $A_-$ to its negative.  \tratsnm
The automorphism has order $1$ precisely when $A_+=A$ and
$A_-=0$. \done
\trats

\subsection{Semisimple idempotents}

As mentioned above, for $a \in A$ and $\lam \in \bbF$, the
$\lam$-eigenspace for the adjoint $\ad_a$ will be qdenoted $A_\lam(a)$
(allowing $A_\lam(a)=0$).

If $A$ is associative and $a$ an idempotent, then $A=A_1(a)\oplus
A_0(a)$.  If $A$ is not associative, then its idempotents can have
adjoint eigenvalues other than $1$ and $0$ and the minimal polynomial
of the adjoint need not be squarefree.  A \alrt{Peirce decomposition}
for $A$ with respect to the \alrt{semisimple idempotent} or
\alrt{axis} $a$ is a decomposition
\[
A = \bigoplus_{\lam \in \Lam} A_\lam(a)\,,
\]
where $\Lam \subseteq \bbF$ is a set containing all eigenvalues for
the adjoint action of $a$ on the algebra $A$.

In particular, Griess and Majorana algebras provide commutative
algebras (typically over $\bbR$) generated by axes whose corresponding
Peirce decompositions have fusion table:
\renewcommand{\arraystretch}{1.75}
\begin{center}
\begin{tabular}{c||ccc|c}
$\star$ & 1  & 0 & $\nfrac{1}{4}$ & $\nfrac{1}{32}$\\
\hline\hline
1 & 1  & $\emptyset$ & $\nfrac{1}{4}$ & $\nfrac{1}{32}$\\
0 & $\emptyset$ & 0 & $\nfrac{1}{4}$ & $\nfrac{1}{32}$\\
$\nfrac{1}{4}$ &
$\nfrac{1}{4}$ & $\nfrac{1}{4}$ & 1,0 & $\nfrac{1}{32}$\\
\hline
$\nfrac{1}{32}$ &
$\nfrac{1}{32}$ & $\nfrac{1}{32}$ & 
$\nfrac{1}{32}$ & 1,0,$\nfrac{1}{4}$\\
\end{tabular}
\end{center}
It is usual \cite[p.\,210]{iv-book} to take $1 \star 0 =0$ in these
rules, but primitivity of the algebra (see \ref{hyp-prim} below)
allows the stronger $1 \star 0 =\emptyset$, which we prefer.

The axes $a$ in these algebras have certain properties that will be of
general interest:

\start{\hyp}{hyp-prim}
\name{Primitivity}
$A_1(a) = \bbF a$.
\trats

\start{\hyp}{hyp-seress}
\name{Seress Condition}
For $\lam \neq 1$, 
\(
A_\lam(a)A_0(a)\subseteq
A_\lam(a) 
\)
\trats

The Seress Condition implies the weaker

\start{\hyp}{hyp-0-subalg}
\name{$0$-subalgebra}
$A_0(a)$ is a subalgebra.
\trats

\start{\hyp}{hyp-zt-grading}
\name{$\zt$-grading}
The eigenvalue set $\Lam$ is the disjoint union of $\Lam_+$ and
$\Lam_-$ with $1 \in \Lam_+$ and such that
\[
A_+ = \bigoplus_{\lam \in \Lam_+} A_\lam(a)
\quad\text{and}\quad
A_- = \bigoplus_{\lam \in \Lam_-} A_\lam(a)
\]
provides a $\zt$-grading of $A$.
\trats

For example, the axes $a$ above have a fusion rule that is a
refinement of the $\zt$-grading
\[
A_+(a) = A_1(a) \oplus A_0(a) \oplus A_{\frac{1}{4}}(a)
\quad\text{and}\quad
A_-(a)= A_{\frac{1}{32}}(a)\,.
\]
Furthermore, when $A_{\frac{1}{32}}(a)=0$ we have
the alternative $\zt$-grading
\[
A_+(a) = A_1(a) \oplus A_0(a)
\quad\text{and}\quad
A_-(a)= A_{\frac{1}{4}}(a)\,.
\]

For the axis $a$ having a $\zt$-grading as in \ref{hyp-zt-grading}, we
define $\tau(a)$ to be the linear transformation of $A$ that acts as
the identity on $A_+(a)$ and negates everything in $A_-(a)$.  By
Proposition \ref{prop-automorphism} these are automorphisms. Those of
order $2$ are the \alrt{Miyamoto involutions} of $A$, since in the
Griess algebra context they were introduced by Miyamoto
\cite{Miyamoto}.  There they were called $\sigma$- or
$\tau$-\alrt{involutions} depending, respectively, upon whether
$A_{\frac{1}{32}}(a)$ was zero or not.  In the Majorana algebra case
these are the \alrt{Majorana involutions} of Ivanov \cite{iv-book}.
The $A_{\frac{1}{32}}(a)=0$ case is related to the \alrt{Atkin-Lehner
  involutions} of Norton \cite[Theorem~3]{No96}.

\subsection{Axial algebras}

An \alrt{axial algebra} is a commutative algebra that is generated by
a set of axes. This has little force unless the axes are provided with
rigid fusion properties, as is the case with the Griess and Majorana
algebras, whose fusion table appears above.

The Griess and Majorana algebra examples are \alrt{primitive axial
  algebras}, in that for each generating axis $a$ we have $A_1(a) =
\bbF a$, as in \ref{hyp-prim}.

Generally, the fusion rule (and table) for the axial algebra $A$
generated by the axis set $\axisset$ has index set $\Lam$, the union
of the eigenvalue sets for the individual $a \in \axisset$, and each
$\mu \star \nu$ satisfies
\[
A_\mu(a) A_\nu(a) \subseteq 
 \bigoplus_{\lam \in \mu\star \nu} A_\lam(a)\,,
\]
for all $a \in \axisset$.

The fusion rule given above for the Griess and Majorana algebras is
actually the particular example $\fusion[V](4,3)$ from the class of
\alrt{Virasoro fusion rules} $\fusion[V](p,q)$, which arise from the
representation theory of rational Virasoro algebras. (See \cite{RS}
for a more detailed discussion.)

Virasoro rules are associated with vertex operator algebras, and the
corresponding Griess algebras come equipped with a form satifying the
Frobenius property:

\start{\hyp}{hyp-frobenius}
\name{Frobenius Property}
There is a bilinear form 
\[
\lla \cdot,\cdot \rra \colon A \times A \goesto \bbF
\]
with $\lla a,a \rra \neq 0$, for all $a \in \axisset$,
that is \alrt{associative:}
\[
\lla ax, b \rra = \lla a, xb \rra
\]
for all $a,b,x \in A$.
\trats

\noindent
For Griess algebras this form is in fact positive definite, but that
supposes characteristic zero and will not appear again until the very
end of the paper.

Frobenius axial algebras have some useful properties.

\start{\prop}{prop-frob}
Let $A$ be a Frobenius axial algebra.
\startnm{a}
\item
The form $\lla \cdot,\cdot \rra$ is symmetric.
\item
For axis $a$ and distinct $\lam,\mu$, the eigenspaces $A_\lam(a)$ and
$A_\mu(a)$ are perpendicular.
\item
If $A$ is primitive, then the radical of the form is the unique
largest ideal of $A$ that does not contain any of the generating axes.
\tratsnm 
\trats

\pf The first two are elementary and can be found in
\cite[Prop.~3.5-6]{RS}.

Each ideal is the direct sum of its eigenspaces for any given
axis. Therefore by \nm{b}, if an ideal of a primitive axial algebra
does not contain any of the generating axes, it is perpendicular to
those axes. In view of associativity of the form, it now follows that
the ideal is perpendicular to any product of the generating axes and
so is in the radical $R$.

On the other hand, for all $a,x\in A$ and $r \in R$
\[
0=\lla ax,r \rra =\lla a,xr \rra \,.
\]
Thus the radical $R$ is itself an ideal and so the maximal ideal
containing none of the axes. \done

\dspace{.5}

The associative algebra $\bigoplus_{i\in I} \bbF a_i$ is a primitive
Frobenius axial algebra generated by the primitive orthogonal
idempotents $\set{a_i}{i \in I}$ and having fusion table
\begin{center}
\begin{tabular}{c||cc}
$\star$ & $1$  & $0$ \\
\hline\hline
$1$ & $1$  & $\emptyset$ \\
$0$ & $\emptyset$ & $0$ 
\end{tabular}
\end{center}
as $a_ia_j=0$ for $i \neq j$.

We can characterize these associative algebras in axial terms.

\start{\prop}{prop-0-1} 
Let $A$ be a primitive axial algebra over the field $\bbF$ that is
generated by the set $\axisset$ of axes.  Then for each $a
\in\axisset$
\[
(\bbF a \oplus A_0(a)) \cap\axisset=  \{a\}
\cup (A_0(a)\cap\axisset)\,.
\]
\trats

\pf 
Let $b \in (\bbF a \oplus A_0(a))\cap\axisset$, 
and set $b=\mu a + n$ for
$\mu \in \bbF$ and $n \in A_0(a)$. Then
\[
\mu a + n = b = b^2 = (\mu a + n)^2 = \mu^2 a +n^2\,,
\]
hence $n=n^2$ and $\mu=\mu^2$. That is, either $\mu=0$ and $n=b \in
A_0(a)$, as desired, or $\mu=1$.

Suppose that $\mu =1$ and $b=a+n$. Therefore
\[
ab=a(a+n)=a^2+an = a \in A_1(b)\,.
\]
But by the primitivity assumption $A_1(b)=\bbF b$ 
contains the unique idempotent $b=a$. \done

\start{\cor}{cor-assoc}
Let $A$ be a primitive axial algebra over the field $\bbF$ that is
generated by the set $\axisset$ of axes having fusion table
\begin{center}
\begin{tabular}{c||cc}
$\star$ & $1$  & $0$ \\
\hline\hline
$1$ & $1$  & $\emptyset$ \\
$0$ & $\emptyset$ & $0$ 
\end{tabular}
\end{center}
\noindent
Then $A = \bigoplus_{a \in \axisset} \bbF a$ is
associative.
\trats

\pf For each $a \in \axisset$, the algebra $A$ is the direct sum of
$\bbF a$ and the subalgebra of $A_0(a)$ generated by all the axes of
$\axisset$ except for $a$. \done

\dspace{1}

In this article we are mainly interested in a case where the algebras
are minimally nonassociative.  For $\eta \in \bbF$ with $0 \neq \eta
\neq 1$, the axial $\bbF$-algebra $A$ is said to have \alrt{Jordan
  type $\eta$} when it is generated by a set of axes with fusion
table:
\begin{center}
\begin{tabular}{c||cc|c}
$\star$ & $1$  & $0$ & $\eta$\\
\hline\hline
$1$ & $1$  & $\emptyset$ & $\eta$\\
$0$ & $\emptyset$ & $0$ & $\eta$\\
\hline
$\eta$ &
$\eta$ & $\eta$ & $1,0$\\
\end{tabular}
\end{center}
In particular, the algebras of Jordan type have a $\zt$-grading and
the Seress Condition.  Indeed for primitive algebras the weaker
$0$-subalgebra condition suffices.

\start{\lem}{lem-subalg-bis}
The primitive axial algebra $A$ has Jordan type $\eta\ (\neq 0,1)$ if
and only if it has the following two properties for each $a$ in its
generating axis set:
\startnm{i}
\item
\name{$0$-subalgebra}
$A_0(a)$ is a subalgebra.
\item
\name{$\zt$-grading}
$A_+(a) =  A_1(a) \oplus A_0(a)$
and
$A_-(a)=A_\eta(a)$
provide a $\zt$-grading of $A$. \done
\tratsnm
\trats

Thus the primitive axial algebras of Jordan type $\eta$ are exactly
the commutative algebras generated by a set of idempotents having the
properties \nm{a}-\nm{d} of the introduction.

\section{Examples}
\label{sec-eg}

\start{}{eg-oA}
{\rm
An axial algebra over $\bbF$ generated by a single axis
is isomorphic to $\bbF$ with the axis $z_0=1$.
This algebra is denoted $\oA$ and only contains
one idempotent, namely its identity element.
}
\trats
\start{}{eg-tB}
{\rm The $2$-dimensional associative algebra $\bbF \oplus \bbF$ with
  $\{b_0,b_1\}$ as $\bbF$-basis and having relations
\[
b_0^2=b_0\,,\ b_1^2=b_1\,,\ b_0b_1=0
\]
is an axial algebra with axes $\{b_0,b_1\}$ and has Jordan type $\eta$
for all $\eta$.  This algebra is denoted $\tB$. Its only idempotents,
other than $b_0$ and $b_1$, are the identity element $b_0+b_1$ and
$0$.  }
\trats
\start{}{eg-tC}
{\rm For $\eta \in \bbF$ let the $\bbF$-algebra $\tC(\eta)$ have basis
  $\{c_0,c_1,c_2\}$ and be subject (for $\{i,j,k\}=\{0,1,2\}$) to
  relations
\[
c_i^2=c_i\,,\ c_ic_j = \nfrac{\eta}{2}(c_i+c_j-c_k).
\]

Then
\begin{align*}
c_i(\eta c_i -c_j -c_k) 
&=\eta c_i^2 -\nfrac{\eta}{2}(c_i + c_j - c_k)
 -\nfrac{\eta}{2}(c_i + c_k - c_j)\\
&=\eta c_i -\nfrac{\eta}{2}c_i 
-\nfrac{\eta}{2}  c_j +\nfrac{\eta}{2} c_k
-\nfrac{\eta}{2}c_i 
-\nfrac{\eta}{2}  c_k +\nfrac{\eta}{2} c_j\\
&=0\,,
\end{align*}
and
\begin{align*}
c_i(c_j -c_k) 
&=\nfrac{\eta}{2}(c_i + c_j - c_k)
 -\nfrac{\eta}{2}(c_i + c_k - c_j)\\
&=\nfrac{\eta}{2}c_i 
+\nfrac{\eta}{2}  c_j -\nfrac{\eta}{2} c_k
-\nfrac{\eta}{2}c_i 
-\nfrac{\eta}{2}  c_k +\nfrac{\eta}{2} c_j\\
&=\eta(c_j-c_k)\,.
\end{align*}

Thus for $\eta \neq 0,1$, the algebra of type $\tC(\eta)$ over $\bbF$
is a $3$-dimensional axial algebra with axis set $\{c_0,c_1,c_2\}$ as
an $\bbF$-basis and generated by $c_0,c_1$.

Comparing this definition with \cite[Table~4]{ipss}, we see that the
algebra $\tC(\nfrac{1}{32})$ over $\bbR$ is the $\tC$ Majorana
algebra, and $\tC(\nfrac{1}{4})$ over $\bbR$ is the $2A$ Majorana
algebra.

The axial algebra $\tC(\eta)$ is, in fact, of Jordan type $\eta$. For
each of the axes $c_i$, the map that fixes $c_i$ and switches $c_j$
and $c_k$ clearly induces an automorphism.  The above calculations
show that this is the Miyamoto involution $\tau(c_i)$, so we have the
needed $\zt$-grading. It remains to prove that the $0$-eigenspace for
$\ad_{c_i}$ is a subalgebra:
\begin{align*}
(\eta c_i -c_j -c_k) &(\eta c_i -c_j -c_k)\\
&=(-c_j -c_k) (\eta c_i -c_j -c_k)\\
&=-\eta c_jc_i -\eta c_kc_i +2c_jc_k
 + c_j +c_k\\
&=-\nfrac{\eta^2}{2}(c_j+c_i-c_k)  -\nfrac{\eta^2}{2}(c_k+c_i-c_j)
\\
&\qquad + 
\eta(c_j+ c_k -c_i)
+c_j +c_k\\
&=-({\eta^2}+{\eta})c_i+(\eta+1)c_j +(\eta+1)c_k\\
&=-(\eta+1)({\eta}c_i-c_j -c_k)\,,
\end{align*}
as desired.

The algebra $\tC(\eta)$ can contain other idempotents, nevertheless
$\bar{\axisset}=\{c_0,c_1,c_2\}$.  It is also a Frobenius algebra, the
symmetric associative bilinear form (unique up to scalar multiple)
being given by
\[
\lla c_i,c_i \rra = 1\,,\
\lla c_i,c_j \rra = \nfrac{\eta}{2}\,.
\]
See Lemma \ref{lem-line-frob} below.

The values $\eta=0,1$ are indeed exceptional.  The algebra $\tC(0)$ is
a copy of the associative algebra $\bbF^3$ of type $\oA^3$ and is not
two-generated.  The algebra $\tC(1)$ has $1$-eigenspaces of dimension
$2$ for its idempotents $c_i$.  It is $\bbF^3$ again, but the
idempotents $c_i$ are now the elements of
$\{(1,1,0),(1,0,1),(0,1,1)\}$.  }
\trats
\start{}{eg-tCs}
{\rm In the special case $\eta=-1$, the algebra $\tC(-1)$ has the
  radical ideal $\bbF(c_0+c_1+c_{2})$.  Let $\tC(-1)^*$ be the
  quotient by this ideal, and set $d_i=c_i +\bbF(c_0+c_1+c_{2})$.
  This is a $2$-dimensional axial algebra spanned by any two of its
  three axes $\{d_0,d_1,d_{2}\,(=-d_0-d_1)\}$ and having
  multiplication given by (for $\{i,j,k\}=\{0,1,2\}$)
\[
d_i^2=d_i\,,\ d_id_j = d_k\,.
\]
As the ideal $\bbF(c_0+c_1+c_{2})$ in $\tC(-1)$ is the radical,
$\tC(-1)^*$ remains a Frobenius algebra.\footnote{Alberto Elduque
  observed that $\tC(-1)^*$ has a natural interpretation as
  para-Hurwitz algebra: depending upon whether or not $\bbF$ contains
  a primitive cube root of unity $\omega$, let $H$ be the Hurwitz
  algebra $\bbF \oplus \bbF$ or $\bbF(\omega)$; then $\tC(-1)^*$ is
  the para-Hurwitz algebra obtained by replacing the usual
  multiplication $\cdot$ for $H$ by $a \star b = (a \cdot b)^J$ where
  $J$ is conjugation in $H$.  Jens K\"oplinger noted that for $\bbF
  \le \bbR$ this para-Hurwitz algebra and its associated form can be
  described nicely in terms of the $A_2$ root lattice thought of a
  subset of the complex plane.  We thank both for their remarks.}  }
\trats
\start{}{eg-jordan}\label{prop-clif-prop}
{\rm A \alrt{Jordan algebra} \cite{jac-jordan} over $\bbF$ (of
  characteristic not two) is a commutative $\bbF$-algebra satisfying
  the identical relation
\[
(a^2 b) a = a^2(ba)\,.
\]

The basic example of a Jordan algebra begins with an associative
algebra $(A,\cdot,+)$ and defines on $A$ a new multiplication (the
\alrt{Jordan product})
\[
a \circ b = \nfrac{1}{2}(ab + ba)\,;
\]
the algebra $(A,\circ,+)$ is then a Jordan algebra, typically denoted
$A^+$.  Those Jordan algebras that arise as subalgebras of $A^+$ for
some associative algebra $A$ are called \alrt{special}.

We are specifically interested in the Clifford algebra
$\Clif{V,\bb{}}$ of the symmetric bilinear form $b$ defined on the
$\bbF$-space $V$.  (See \cite[p.75]{jac-jordan} for discussion of
Clifford algebras.)  $\Clif{V,\bb{}}$ is the associative
$\bbF$-algebra that results from factoring out of the tensor algebra
$\operatorname{T}(V)$ the ideal generated by all elements $v \otimes w
+ w \otimes v - b(v,w) 1$ for $v,w \in V$.  As $v \otimes v -
\nfrac{1}{2}b(v,v)$ is in the ideal, the Clifford algebra is spanned
by the various monomials $\prod_{i=1}^k v_{\iota(i)}$, for nonnegative
$k$, where $\iota$ is an injection of $1, \dots, k$ into the index set
$I$ for the basis $\set{v_i}{i \in I}$ of $V$.  The Clifford algebra
admits the canonical reversal involution $J$ determined by
\[
\sym\colon\prod_{i=1}^k v_{\iota(i)} \mapsto \prod_{i=k}^1 v_{\iota(i)}
\]
for all $k$ and all $k$-subsets $\set{v_{\iota(i)}}{1 \le i \le k}$.
Its fixed points form $\SClif{V,\bb{}}$, a Jordan subalgebra of
$\Clif{V,\bb{}}^+$.
\label{eg-cliff}

The special Jordan algebra $\SClif{V,\bb{}}$ of $\sym$-symmetric
Clifford elements in $\Clif{V,\bb{}}$ contains the subset
$\Vp{V}{\bb{}}= \bbF 1 \oplus V$ which is in fact a (special) Jordan
subalgebra itself because
\begin{align*}
(\alpha 1 + v) \circ (\beta 1 + w) 
=& \nfrac{1}{2}\left((\alpha 1+v)(\beta 1+w)
+(\beta 1+w)(\alpha 1+v)\right)\\
=& \nfrac{1}{2}\left(\alpha\beta 1+\alpha w +\beta v + vw
+\beta\alpha 1+\beta v + \alpha w + wv\right)\\
=& \nfrac{1}{2}\left(2\alpha\beta 1+2\alpha w 
+2\beta v + vw + wv\right)\\
=& \nfrac{1}{2}\left((2\alpha\beta+b(v,w))1 
+2\alpha w +2\beta v\right)\\
=& (\alpha\beta+\nfrac{1}{2}b(v,w))1 +\alpha w +\beta v\,.
\end{align*}

We collect some properties of the Jordan algebra $V^\sym$:

{\em
\startnm{a}
\item The nonidentity idempotents of $V^\sym$ are exactly the elements
  $\frac{1}{2} + v$ for $v \in V$ with $b(v,v)=\frac{1}{2}$.
\item If $e^+=\frac{1}{2} + v$ is a nonidentity idempotent, then so is
  $e^-=\frac{1}{2} - v$.  Furthermore for $\eps = \pm$
\[
V^\sym_1(e^\eps)=\bbF e^\eps\,,\
V^\sym_0(e^\eps)=\bbF e^{-\eps}\,,
\]
and
\[
V^\sym_{\frac{1}{2}}(e^+)
=V^\sym_{\frac{1}{2}}(e^-)
=v^\perp
\,.
\]
In particular the idempotent $e^\eps$ is semisimple with eigenvalues
$\{1,0,\nfrac{1}{2}\}$ and $1$- and $0$-eigenspaces of dimension $1$.
\item
For $e^\eps=\frac{1}{2} + \eps v\ (\eps=\pm)$ a nonidentity
idempotent, let $t$ with $t(1)=1$ be the extension to $V^\sym$ of the
negative of the orthogonal reflection on $V$ with center $v$. Then $t$
is the Miyamoto involution $\tau(e^+)=\tau(e^-)$ on $V^\sym$. In
particular, $V^\sym$ is $\zt$-graded by $e^\eps$ with $\Lam_+=\{1,0\}$
and $\Lam_-=\{\nfrac{1}{2}\}$.  
\tratsnm 
}

{\rm \pf We identify $\bbF$ and $\bbF1$. Also we no longer use $\circ$
  to indicate multiplication in $V^\sym$, which is thus given by
\begin{align*}
(\alpha  + v)(\beta  + w) 
&= (\alpha\beta+\nfrac{1}{2}b(v,w)) +\alpha w +\beta v\,.
\end{align*}
\startnm{a}
\item
We have 
\[
(\alpha + v)^2 = (\alpha^2 + \nfrac{1}{2} b(v,v)) + 2\alpha v\,.
\]
For idempotent $\alpha +v$ this forces $\alpha=\frac{1}{2}$ and then
$b(v,v) = \nfrac{1}{2}$.
\item
As $b(v,v)=b(-v,-v)$, $e^+$ is an idempotent if and only if $e^-$ is,
by the previous part.  Next
\[
e^\eps(\beta  + w) =
(\nfrac{1}{2}  + \eps v)(\beta  + w) 
= (\nfrac{\beta}{2}+\nfrac{\eps}{2}b(v,w)) +\nfrac{1}{2}w 
+\beta\eps v\,.
\]
For $\beta + w$ to be in $V^\sym_1(e^+)$ we must have
$\nfrac{1}{2}w +\beta\eps v=w$; that is, $w=2\beta\eps v$.
Therefore 
\[
V^\sym_1(e^\eps) = \{\beta + 2 \beta\eps v\}
=\{2\beta(\nfrac{1}{2}+\eps v) \} =\bbF e^\eps\,.
\]
Similarly if $\beta + w \in V^\sym_0(e^\eps)$ then
$\nfrac{1}{2}w +\beta\eps v= 0$ and $w=-2\beta\eps v$,
leading to $V^\sym_0(e^\eps)=\bbF e^{-\eps}$.

Let $W=v^{\perp}$, the subspace of $V$ consisting of elements that are
$b$-perpendicular to $v$. Then
\[
V^\sym= \bbF e^+\oplus  \bbF e^- \oplus W = \bbF 1 
\oplus \bbF v \oplus W\,.
\]
Then for $w \in W$
\[
e^\eps w = (\nfrac{0}{2} + \nfrac{\eps}{2}b(v,w))
+ \nfrac{1}{2}w
+ 0\cdot \eps v =\nfrac{1}{2}w\,,
\]
so that $W=V^\sym_{\frac{1}{2}}(e^\eps)$
and
\[
V^\sym= V^\sym_{1}(e^\eps) \oplus
V^\sym_{0}(e^\eps) \oplus
V^\sym_{\frac{1}{2}}(e^\eps)\,.
\]

\item The negative reflection in $v$ is an isometry of the form $b$ on
  $V$; so its extension $t$ to $V^\sym$ is an algebra automorphism,
  which by the previous part is equal to both Miyamoto involutions
  $\tau(e^+)$ and $\tau(e^-)$. The existence of this automorphism is
  equivalent to $\zt$-grading.  
\done 
\tratsnm

Consider the symmetric Jordan Clifford algebra
$\SClif{\bbF^2,\bb{\delta}}$ on the vector space $\bbF^2$ with basis
$v_0,v_1$ and symmetric bilinear form given by
$\bb{\delta}(v_i,v_i)=2$, $\bb{\delta}(v_0,v_1)=\delta$ for some
constant $\delta \in \bbF$. (Note that different choices of $\delta$
may give isometric spaces and so isomorphic algebras.)

As $V$ has dimension $2$, its Clifford algebra is spanned by $1, v_0,
v_1, v_0v_1$; so the Jordan algebra $\SClif{\bbF^2,\bb{\delta}}$ is
equal to its subalgebra $V^\sym$ with basis $1, v_0, v_1$.

This is an axial algebra with generating axes
\[
e_0^+=\nfrac{1}{2}(1+v_0)\,,\
e_0^-=\nfrac{1}{2}(1-v_0)\,,\
e_1^+=\nfrac{1}{2}(1+v_1)\,,\
e_1^-=\nfrac{1}{2}(1-v_1)\,,\
\]
or, in fact, any three of these. Furthermore
\begin{align*}
e_0^+ e_1^+ &= \nfrac{1}{2}(1+v_0)\nfrac{1}{2}(1+v_1)\\
&= \nfrac{1}{4}(1+v_0+v_1+v_0v_1)\\
&= \nfrac{1}{4}(1+v_0+v_1+\nfrac{1}{2}b(v_0,v_1))\\
&= \nfrac{1}{4}(1+v_0+v_1+\nfrac{\delta}{2})\\
&= \nfrac{1}{4}(1+v_0+1+v_1+(\nfrac{\delta}{2}-1))\\
&= \nfrac{1}{2}e_0^++\nfrac{1}{2}e_1^+ +\nfrac{1}{8}(\delta-2)1\,;
\end{align*}
so $V^\sym$ is generated by $e_0^+$ and $e_1^+$ except when
$\delta=2$.

Every Jordan algebra generated by idempotents is a Frobenius algebra
\cite{jac-jordan}.  Here that is easy to check for the form given by
\begin{center}
\begin{tabular}{c|ccc}
$\lla \cdot, \cdot \rra$   
& $1$ & $v_0$ & $v_{1}$\\
\hline
$1$ & $2$ & $0$ &  $0$\\
$v_0$ & $0$ &$2$ &$\delta$\\
$v_{1}$ & $0$ & $\delta$& $2$\\
\end{tabular}\quad
\end{center}
That is, the form $b_\delta$ on $V$ is extended to the full algebra
$\bbF 1 \oplus V$ by setting $\lla 1,1\rra =2$ and $V = 1^\perp$.

\dspace{.5}

For $\delta=2$, the idempotents $e_0^+$ and $e_1^+$ generate a
$2$-dimensional Jordan subalgebra $\SClif[0]{\bbF^2,\bb{2}}$ that is
of type $\oA$ when factored by the $1$-dimensional ideal spanned by
$e_1^+ - e_0^+=\nfrac{1}{2}(v_1-v_0)$.

In any event we have the multiplication table:
\begin{center}
\begin{tabular}{c||c|c|c}
$\centerdot $   
& $e_0^+$ & $e_1^+$  & $s$\\
\hline\hline
$e_0^+$  &$e_0^+$  & $\onehalf e_0^+ + \onehalf e_1^+ +s$ 
&  $pe_0^+$ \\
$e_1^+$  & $\onehalf e_0^+ + \onehalf e_1^+ + s$ &$e_1^+$  
& $p e_1^+$ \\
$s$ & $pe_0^+$  & $pe_1^+$ & $ps$\\
\end{tabular}\quad
\end{center}
for $p=\nfrac{1}{8}(\delta -2) \in \bbF$ and
$s=\nfrac{1}{8}(\delta -2)1 \in V^J$.

\dspace{.5}

Isometric forms $b_\delta,\,b_{\delta'}$ produce isomorphic Clifford
algebras, hence isomorphic symmetric Jordan Clifford algebras
$\SClif{\bbF^2,\bb{*}}$. For instance, the algebras
$\SClif{\bbF^2,\bb{\delta}}$ and $\SClif{\bbF^2,\bb{-\delta}}$ are
isomorphic as they correspond, respectively, to the bases $v_0$, $v_1$
and $v_0'=v_0$, $v_1'=-v_1$.  This amounts to replacing the pair of
idempotents $e_0^+,e_1^+$ by the pair $(e_0')^+=e_0^+,
(e_1')^+=e_1^-$. The above multiplication table makes it clear that no
automorphism of $\SClif{\bbF^2,\bb{\delta}}$ takes the pair of axes
$\axisset=\{e_0^+,e_1^+\}$ to the pair $\axisset'=\{(e_0')^+=e_0^+,
(e_1')^+=e_1^-\}$.

Indeed in the exceptional degenerate case $\delta=2$, the axis set
$\axisset$ does not even generate $\SClif{\bbF^2,\bb{\delta}}$ (as
noted above) while $\axisset'$ does.  }
}\trats

\start{}{eg-cliff00}
{\rm In the previous example, when $\delta=2$ the element $s$ equals
  $0$ since the scalar $p=\onehalf(\delta-2)$ is $0$.  But we can
  adapt the corresponding multiplication table to include that
  case. The algebra $\SClif[00]{\bbF^2,\bb{2}}$ is the vector space
  $\bbF e_0' \oplus \bbF e_1' \oplus \bbF s'$ subject to the
  multiplication table:
\begin{center}
\begin{tabular}{c||c|c|c}
$\centerdot $   
& $e_0'$ & $e_1'$  & $s'$\\
\hline\hline
$e_0'$  &$e_0'$  & $\onehalf e_0' + \onehalf e_1' +s'$ &  $0$ \\
$e_1'$  & $\onehalf e_0' + \onehalf e_1' + s'$ &$e_1'$  & $0$ \\
$s'$ & $0$  & $0$ & $0$\\
\end{tabular}\quad
\end{center}
This is a Jordan algebra whose ideal $\bbF s'$ is the $0$-eigenspace
for both of the generating idempotents $e_0'$ and $e_1'$.  The
quotient by this ideal is a copy of $\SClif[0]{\bbF^2,\bb{2}}$.  }
\trats

\section{Two-generated axial algebras of Jordan type}
\label{sec-poc}

Let $A$ be a primitive axial algebra of Jordan type $\eta$ over $\bbF$
that is generated by the two axes $a$ and $b$.

For $x \in A$ and $p \in \axisset=\{a,b\}$, write
\[
x = \phi_p(x)p + \alpha_p(x) + \gamma_p(x)
\]
for $\phi_p(x) \in \bbF$, $\alpha_p(x) \in A_0(p)$, and $\gamma_p(x)
\in A_\eta(p)$.

Recall that the Miyamoto involutions $\tau(p)$ are automorphisms that
act via
\[
x^{\tau(p)} = \phi_p(x)p + \alpha_p(x) - \gamma_p(x)\,.
\]

\start{\lem}{lem-ab}
\startnm{a}
\item $\gamma_b(a)=\onehalf(a - a^{\tau(b)})$ and
      $\gamma_a(b)=\onehalf(b - b^{\tau(a)})$.
\item $\alpha_b(a)=-\phi_b(a) b + \onehalf(a + a^{\tau(b)})$ and
      $\alpha_a(b)=-\phi_a(b) a + \onehalf(b + b^{\tau(a)})$.
\item $ab=\phi_b(a)b + \eta \gamma_b(a)$ and
      $ab=\phi_a(b)a + \eta \gamma_a(b)$.
\item Set $\sigma=ab- \eta a - \eta b$. Then
\[
\sigma=
((1-\eta)\phi_b(a)-\eta)b -\eta \alpha_b(a)
=((1-\eta)\phi_a(b)-\eta)a -\eta \alpha_a(b)
\,.
\]
In particular $\sigma$ is fixed by both $\tau(a)$ and $\tau(b)$.
\tratsnm
\trats

\pf
By symmetry, we need only prove one equality from each part.

\noindent
\nm{a}
We have
\[
a = \phi_b(a)b + \alpha_b(a) + \gamma_b(a)
\]
hence
\[
a^{\tau(b)} = \phi_b(a)b + \alpha_b(a) - \gamma_b(a)\,,
\]
so $a-a^{\tau(b)} = 2 \gamma_b(a)$.

\noindent
\nm{b}
Next
\begin{align*}
a + a^{\tau(b)}
&=
\left(\phi_b(a)b + \alpha_b(a) + \gamma_b(a)\right)
+\left(\phi_b(a)b + \alpha_b(a) - \gamma_b(a)\right)\\
&=2(\phi_b(a)b + \alpha_b(a))\,.
\end{align*}

\noindent
\nm{c}
Furthermore
\[
ab = (\phi_b(a)b + \alpha_b(a) + \gamma_b(a))b
=\phi_b(a)b + \gamma_b(a)b=\phi_b(a)b + \eta\gamma_b(a)\,.
\]

\noindent
\nm{d}
Finally
\begin{align*}
\sigma &= ab- \eta a - \eta b\\
&=(\phi_b(a)b + \eta \gamma_b(a)) 
-\eta (\phi_b(a)b + \alpha_b(a) + \gamma_b(a))
- \eta b\\
&=
((1-\eta)\phi_b(a)-\eta)b -\eta \alpha_b(a) 
\,,
\end{align*}
which is in $A_1(b)+A_0(b)$ and so is fixed by $\tau(b)$.
\done

\start{\lem}{lem-a-sigma}
\startnm{a}
\item $a\sigma = ((1-\eta)\phi_a(b)-\eta)a$.
\item $b\sigma = ((1-\eta)\phi_b(a)-\eta)b$.
\tratsnm
\trats

\pf We need only consider \nm{b}.
From Lemma \ref{lem-ab}
\[
b\sigma
= b \left(((1-\eta)\phi_b(a)-\eta)b -\eta \alpha_b(a)\right)
= ((1-\eta)\phi_b(a)-\eta)b
\,,
\]
as claimed.  \done

\start{\lem}{lem-seress}\name{Seress Lemma}\ \
For $p \in \axisset$, $x \in A$, and $y \in A_1(p)+A_0(p)$
\[
p(xy)=(px)y\,.
\]
\trats

\pf
If $z=\alpha p \in A_1(p)$, then 
\[
p(xz)=p(x (\alpha p))=(\alpha p)(xp)
=(xp)(\alpha p)=(px)z\,;
\]
so by linearity we may assume $y \in A_0(p)$.

Also by linearity we may assume that $x \in A_\lam(p)$.  In particular
if $x = \beta p \in A_1(p)$, then
\[
p(xy)=p(\beta p y)=0=
\beta py =(\beta pp)y=(px)y\,.
\]
Finally if $x \in A_\lam(p)$ for $\lam \neq 1$, then by the fusion
rules $xy \in A_\lam$ also; so
\[
p(xy)=\lam(xy)=(\lam x)y =(px)y\,. \quad \done
\]

\start{\lem}{lem-phi}
$(ab)\sigma=\pi ab$ for $\pi = (1-\eta)\phi-\eta\in \bbF$ with
$\phi=\phi_a(b)=\phi_b(a)$.
\trats

\pf If $ab=0$, then certainly $(ab)\sigma=0=\pi ab$.  In that case $b
\in A_0(a)$ and $a \in A_0(b)$ so that also $\phi_a(b)=0=\phi_b(a)$.
Thus we may assume that $ab \neq 0$.

As $\sigma$ is fixed by $\tau(a)$ and $\tau(b)$, we have
\[
\sigma \in (A_1(a)+A_0(a)) \cap (A_1(b) + A_0(b))\,.
\]
Therefore by the previous lemma
\[
b(a\sigma)=(ba)\sigma=(ab)\sigma=a(b\sigma)\,,
\]
and so by Lemma \ref{lem-a-sigma}
\begin{align*}
((1-\eta)\phi_a(b)-\eta)ab&=
b\left(((1-\eta)\phi_a(b)-\eta)a\right)\\
&=b(a\sigma)=a(b\sigma)\\
&=a\left(((1-\eta)\phi_b(a)-\eta)b\right)\\
&=
((1-\eta)\phi_b(a)-\eta)ab\,.
\end{align*}
As $ab \neq 0$
\[
(1-\eta)\phi_a(b)-\eta=
(1-\eta)\phi_b(a)-\eta=\pi
\]
and especially $\phi_a(b)=\phi_b(a)$ since $\eta\neq 1$.
\done

\start{\lem}{lem-sigma-sigma}
$\sigma^2=\pi \sigma$.
\trats

\pf
By Lemmas \ref{lem-a-sigma} and \ref{lem-phi}
\[
\sigma^2=(ab-\eta a -\eta b)\sigma
=\pi ab - \eta \pi a -\eta \pi b = \pi \sigma\,.\quad \done
\]

\dspace{1}

\start{\prop}{prop-dim3}
$A = \bbF a + \bbF b + \bbF \sigma$
with multiplication table:
\begin{center}
\begin{tabular}{c||c|c|c}
$\centerdot $   
& $a$ & $b$ & $\sigma$\\
\hline\hline
$a$ & $a$ & $\eta a+ \eta b + \sigma$ &  $\pi a$\\
$b$ & $\eta a+\eta b + \sigma$ &$b$ &$\pi b$\\
$\sigma$ & $\pi a$ & $\pi b$& $\pi \sigma$\\
\end{tabular}\quad
\end{center}
where $\pi = (1-\eta)\phi-\eta$.
\trats

\pf The multiplication table is immediate from the previous lemmas and
our definitions of $a$, $b$, and $\sigma$. As the span $\bbF a + \bbF
b + \bbF \sigma$ contains the generators of $A$ and contains the
pairwise products of its three spanning elements, it is all of $A$.
\done

\dspace{.5}

We are not claiming in this proposition that $a$, $b$, and $\sigma$
are linearly independent. Indeed in algebras of type $\oA$ and $\tB$
this is certainly not the case.

\start{\thm}{thm-b}
Choose parameters $\eta\,(\neq 0,1)$ and $\phi$ in $\bbF$.  Let $B =
B(\eta,\phi)= \bbF c \oplus \bbF d \oplus \bbF \rho$ be the
commutative $\bbF$-algebra with multiplication table:
\begin{center}
\begin{tabular}{c||c|c|c}
$\centerdot $   
& $c$ & $d$ & $\rho$\\
\hline\hline
$c$ & $c$ & $\eta c+ \eta d + \rho$ &  $\pi c$\\
$d$ & $\eta c+ \eta d + \rho$ &$d$ &$\pi d$\\
$\rho$ & $\pi c$ & $\pi d$& $\pi \rho$\\
\end{tabular}\quad
\end{center}
where
$\pi = (1-\eta)\phi-\eta$. 
\startnm{a}
\item
$B(\eta,\phi)$ is a primitive axial algebra generated by axes
  $\calB=\{c,d\}=\{p,q\}$ with
\[
B=B_1(p) \oplus B_0(p) \oplus B_\eta(p)
\]
for $B_1(p)=\bbF p$, $B_0(p)=\bbF(\pi p - \rho)$, and
$B_\eta(p)=\bbF((\eta-\phi)p + \eta q + \rho)$.
\item
$B(\eta,\phi)$ satisfies the Seress condition:
\[
B_0(c)B_\lam(c)\subseteq B_\lam(c)\ \ \text{for}\ \lam \neq 1.
\]
\item
For the generating set of axes $\calB$ we have the fusion table:
\begin{center}
\begin{tabular}{c||cc|c}
$\star$ & $1$  & $0$ & $\eta$\\
\hline\hline
$1$ & $1$  & $\emptyset$ & $\eta$\\
$0$ & $\emptyset$ & $0$ & $\eta$\\
\hline
$\eta$ &
$\eta$ & $\eta$ & $1,0,\eta$\\
\end{tabular}
\end{center}
\item
$B(\eta,\phi)$ is simple except in the following cases: \startnm{i}
\item $\phi=\nfrac{\eta}{1-\eta}$
and $\pi=0$
where $B_0(c)=B_0(d)=\bbF \rho$
is an ideal;
\item $\phi=0$ and $\pi=-\eta$
where $B_\eta(c)=B_\eta(d)
=\bbF(\eta c + \eta d + \rho)$
is an ideal
as are 
\[
B_0(c)\oplus B_\eta(c)=\bbF d \oplus \bbF(\eta c + \eta d + \rho)
=B_1(d)+B_\eta(d)
\]
and 
\[
B_0(d)\oplus B_\eta(d)=\bbF c \oplus
\bbF(\eta c + \eta d + \rho)=B_1(c)\oplus B_\eta(c)\,;
\]
\item $\phi=1$ and
$\pi=1-2 \eta$ where
\[
B_0(c)\oplus B_\eta(c)=B_0(d)\oplus B_\eta(d)
\]
is an ideal.
\tratsnm
\tratsnm
\trats

\pf \nm{a} Certainly $p \in B_1(p)$.  Also
\[
p(\pi p - \rho)= \pi p^2 - p\rho=\pi p - \pi p = 0,
\]
so $\pi p - \rho \in B_0(p)$. Next
\begin{align*}
p\left((\eta-\phi)p + \eta q + \rho\right) 
&=(\eta-\phi)p^2 + \eta pq + p\rho\\
&=(\eta-\phi)p + \eta (\eta p + \eta q + \rho) + \pi p\\
&=(\eta-\phi+\eta^2 +\pi)p 
+  \eta^2 q +\eta\rho\\
&=(\eta-\phi+\eta^2 + \phi-\eta\phi - \eta)p 
+  \eta^2 q +\eta\rho\\
&=\eta\left((\eta-\phi)p 
+  \eta q +\rho\right)\,,
\end{align*}
hence $(\eta-\phi)p + \eta q + \rho \in B_\eta(p)$.  As $\eta \neq 0$,
the three vectors $p$, $\pi p - \rho$, and $(\eta-\phi)p + \eta q +
\rho$ are linearly independent in $B$ of dimension $3$.

\dspace{.2}

\nm{b} The Seress Condition speaks to the entries of the fusion table
corresponding to $B_0(p)B_\lam(p)$ with $\lam \in \{0,\eta\}$.  Here
$B_0(p)$ is spanned by $\pi p - \rho$.  But each $B_\lam(p)$ is an
eigenspace for $p$ by definition, while $\rho$ acts on $B$ as scalar
multiplication by $\pi$.  Therefore in all cases
$B_0(p)B_\lam(p)\subseteq B_\lam(p)$.

\nm{c} The fusion table summarizes parts of \nm{a} and \nm{b}.

\nm{d} As the adjoint eigenvalues $1$, $0$, and $\eta$ are distinct,
any ideal is a direct sum of certain of the $1$-dimensional
eigenspaces $B_1(c)$, $B_0(c)$, and $B_\eta(c)$ and equally well of
$B_1(d)$, $B_0(d)$, and $B_\eta(d)$.

No ideal of dimension $1$ can contain $c$ or $d$, since the quotient
of dimension $2$ would be generated by a single idempotent. Therefore
a $1$-dimensional ideal is one of $B_0(c)$ or $B_\eta(c)$ that is
simultaneously equal to one of $B_0(d)$ or $B_\eta(d)$.

First consider an ideal $B_0(c)=\bbF(\pi c - \rho)$.  It must also
contain
\begin{align*}
(\pi c - \rho)d &= \pi dc - d \rho\\
&=\pi(\eta c + \eta d + \rho) - \pi d\\
&=\pi(\eta c + (\eta-1) d + \rho)
\,,
\end{align*}
As $\eta\neq 1$, the element $\pi(\eta c + (\eta-1) d + \rho)$ is a
scalar multiple of $\pi c - \rho$ in $B(\eta,\phi)$ if and only if
$\pi=0$ hence $\phi=\nfrac{\eta}{1-\eta}$.  The corresponding ideal of
dimension $1$ is $\bbF \rho$, which is $B_0(d)$ as well.  That is,
$B_0(c)$ is an ideal if and only if $B_0(d)$ is an ideal, in which
case both are $\bbF \rho$.

By the above, $B_\eta(c)$ is an ideal if and only if it equals
$B_\eta(d)$. The generators are, respectively, $(\eta-\phi)c + \eta d
+ \rho$ and $\eta c +(\eta-\phi)d + \rho$, equal precisely when
$\phi=0$.  The corresponding ideal of dimension $1$ is then
\[
B_\eta(c)=
\bbF(\eta c + \eta d + \rho)=B_\eta(d)\,.
\]

An ideal of dimension $2$ in $B(\eta,\phi)$ not containing $c$ must be
\[
I(\eta,\phi)=B_0(c)\oplus B_\eta(c)=\bbF(\pi c - \rho) \oplus
\bbF((\eta-\phi)c + \eta d + \rho)\,.
\]
The element $d$ must map to an idempotent $\bar{d}$ in
$\bar{B}=B(\eta,\phi)/I(\eta,\phi)$, which is isomorphic to axial
algebra $\bbF$ of type $\oA$. Therefore $\bar{d}$ is either $\bar{c}$,
so that $c-d \in I(\eta,\phi)$, or $\bar{0}$, which is to say $d\in
I(\eta,\phi)$.  On the other hand,
\[
d = \phi\, c + \eta^{-1}\left((\pi c - \rho) +
((\eta-\phi)c + \eta d + \rho)\right) \in \phi\, c + I(\eta,\phi)\,;
\]
so the two cases lead, respectively, to $\phi = 1$ and $\phi=0$.

If $\phi=1$, then $\pi=1-2\eta$.  Here the ideal $I(\eta,\phi)$ is
\begin{align*}
B_0(c)\oplus B_\eta(c)&=
\bbF(\pi c - \rho) \oplus
\bbF((\eta-1)c + \eta d + \rho)\\
&=\bbF((1-2\eta) c - \rho) \oplus
\bbF((\eta-1)c + \eta d + \rho)\\
&=\bbF (c-d) \oplus \bbF(-\pi d + \rho)\\
&=B_0(d)\oplus B_\eta(d)\,.
\end{align*}

Next suppose $\phi=0$. A $2$-dimensional ideal that contains $d$ but
not $c$ is thus
\[
B_0(c)\oplus B_\eta(c)=\bbF d \oplus
\bbF(\eta c + \eta d + \rho)=B_1(d)+B_\eta(d)\,.
\]
while by symmetry a $2$-dimensional ideal containing $c$ but not $d$
is
\[
B_0(d)\oplus B_\eta(d)=\bbF c \oplus
\bbF(\eta c + \eta d + \rho)=B_1(c)+B_\eta(c)\,. 
\quad\done
\]

\start{\prop}{prop-image-jordan}
Let $\bar{B}$ be a quotient of $B=B(\eta,\phi)$.  Then $\bar{B}$ is an
axial algebra of Jordan type $\eta$ if and only if we have one of:
\startnm{1}
\item
$\bar{B}$ is associative and isomorphic to $\bbF$ of type $\oA$ or
  $\bbF\oplus\bbF$ of type $\tB$;
\item
$\phi = \nfrac{\eta}{2}$;
\item
$\eta=\onehalf$.
\tratsnm
In all these cases $\bar{B}$ is spanned by $\bar{c}$, $\bar{d}$, and
either of $\bar{c}^{\tau(\bar{d})}$ or $\bar{d}^{\tau(\bar{c})}$.
\trats

\pf By the previous theorem, $\bar{B}$ is of Jordan type $\eta$
precisely when
\[
\bar{B}_\eta(\bar{p})\bar{B}_\eta(\bar{p}) 
\subseteq \bar{B}_1(\bar{p}) \oplus \bar{B}_0(\bar{p})\,.\,
\]
where $\{\bar{p},\bar{q}\}=\{\bar{c},\bar{d}\}$.

In $B=B(\eta,\phi)$ we have $B_\eta(p)=\bbF((\eta-\phi)p + \eta q +
\rho)$.  Thus $\bar{B}_\eta(\bar{p})$ is spanned by
$(\eta-\phi)\bar{p} + \eta \bar{q} + \bar{\rho}$.  To check fusion
containment we calculate
\begin{align*}
((\eta&-\phi)\bar{p} + \eta \bar{q} + \bar{\rho})((\eta-\phi)\bar{p} 
+ \eta \bar{q} + \bar{\rho})\\
&=(\eta-\phi)^2\bar{p}^2+\eta^2 \bar{q}^2 + \bar{\rho}^2
+2\eta \bar{q} \bar{\rho} + 2 (\eta-\phi) \bar{p} \bar{\rho} 
+ 2(\eta-\phi)\eta \bar{p}\bar{q}\\
&=(\eta-\phi)^2\bar{p}+\eta^2 \bar{q} + \pi\bar{\rho}
+2\pi\eta \bar{q}  + 2\pi (\eta-\phi) \bar{p}  
+ 2(\eta-\phi)\eta (\eta \bar{p} + \eta \bar{q} + \bar{\rho})\\
&=((\eta-\phi)^2+ 2\pi (\eta-\phi)+ 2(\eta-\phi)\eta^2)\bar{p} 
+( \eta^2+2\pi\eta+2(\eta-\phi)\eta^2)\bar{q}\\
&\qquad+( \pi+2(\eta-\phi)\eta)\bar{\rho}\,.
\end{align*}
Therefore $\bar{B}$ is an axial algebra of Jordan type $\eta$ if and
only if the subalgebra $\bar{B}_1(\bar{p}) \oplus \bar{B}_0(\bar{p})
=\bbF \bar{p} \oplus \bbF(\pi \bar{p} - \bar{\rho}) =\bbF \bar{p}
\oplus \bbF\bar{\rho}$ contains the element
\begin{align*}
(\eta^2+2\pi\eta+2(\eta-\phi)\eta^2)\bar{q} 
&=(\eta^2+2((1-\eta)\phi-\eta)\eta+2(\eta-\phi)\eta^2)\bar{q}\\
&=(\eta^2+2\phi\eta-2\phi\eta^2-2\eta^2+2\eta^3-2\phi\eta^2)\bar{q}\\
&=\eta(2\eta^2+(-1-4\phi)\eta+2\phi)\bar{q}\\
&=\eta(2\eta-1)(\eta-2\phi)\bar{q}\,.
\end{align*}

As $\bbF$ and $\bbF \oplus \bbF$ have Jordan type $\eta$ for all
$\eta$, this immediately gives us the converse part of the
proposition, the claim about spanning following from Lemma
\ref{lem-ab}(c).

Now assume that $\bar{B}$ does have Jordan type $\eta$ but $\phi \neq
\nfrac{\eta}{2}$ and $\eta \neq \onehalf$.  Then
$\eta(\eta-2\phi)(2\eta-1)$ is nonzero in $\bbF$, and so
\[
\bar{d} \in \bbF\bar{c} \oplus B_0(\bar{c})=\bar{B}
\quad\text{and}\quad
\bar{c} \in \bbF\bar{d} \oplus B_0(\bar{d})=\bar{B}\,.
\]
By Corollary \ref{cor-assoc} the axial algebra generated by $\bar{c}$
and $\bar{d}$ is associative and isomorphic to $\bbF$ or to $\bbF
\oplus \bbF$, as desired.  \done

\dspace{1}

\noindent
{\sc Proof of Theorem \ref{thm-2gen-jordan-field}.}

\noindent
By Propositions \ref{prop-dim3} and \ref{prop-image-jordan} the
algebras $A$ of the theorem are $\bbF$, $\bbF \oplus \bbF$, and the
quotients of $B(\eta,\nfrac{\eta}{2})$ and $B(\onehalf,\phi)$.  We
claim:
\begin{quote}
\nm{i} $B(\eta,\nfrac{\eta}{2})$ is isomorphic to $\tC(\eta)$.

\nm{ii} $B(\onehalf,\phi)$ is isomorphic to
$\SClif{\bbF^2,\bb{\delta}}$ for $\phi \neq 1$, $\delta=4\phi - 2 \neq
2$, and $B(\onehalf,1)$ is isomorphic to $\SClif[00]{\bbF^2,\bb{2}}$.
\end{quote}

\nm{i} As $B(\eta,\nfrac{\eta}{2})$ and $\tC(\eta)$ both have
dimension $3$, it is enough to note that the axial algebra $\tC(\eta)$
has Jordan type $\eta$. But this was shown in \ref{eg-tC}.

\dspace{.5}

\nm{ii} By \ref{eg-cliff} and \ref{eg-cliff00} the $3$-dimensional
algebras $\SClif{\bbF^2,\bb{\delta}}$ ($\delta \neq 2$) and
$\SClif[00]{\bbF^2,\bb{2}}$ ($\delta = 2$) are quotients of the
$3$-dimensional algebras $B(\onehalf,\phi)$ for $\pi=
\nfrac{1}{8}(\delta-2)$. As $\pi=(1-\eta)\phi-\eta =\onehalf \phi -
\onehalf$, this gives $\delta= 4 \phi -2$.

\dspace{.5}

The proper quotients of $B(\eta,\phi)$ are detailed in Theorem
\ref{thm-b}. A quotient of dimension $1$ must have type $\oA$ and need
not be discussed further. A quotient of dimension $2$ with $\phi=0$
has $\bar{B}_\eta(p)=0$ and so is of type $\tB$. Therefore we only
need consider quotients of dimension $2$ with
$\phi=\nfrac{\eta}{1-\eta}$.  This leads to $\bar{B}(-1,-\onehalf)$,
which is isomorphic to $\tC(-1)^*$ by \ref{eg-tC}, and
$\bar{B}(\onehalf,1)$, which is isomorphic to
$\SClif[0]{\bbF^2,\bb{2}}$ by \ref{eg-cliff00}.  \done

\start{\rmks}{rmk-dups}
\nm{1} For an algebra of dimension $3$ to appear under both \nm{i} and
\nm{ii} we must have $\eta=\onehalf$ and $\phi=\frac{1}{4}$ so that
$\delta=4\phi-2=-1$.  Thus the only dimension $3$ algebra to occur in
both is $\tC(\onehalf)$ which is $\SClif{\bbF^2,\bb{-1}}$ when the
characteristic is not three and $\SClif[00]{\bbF^2,\bb{2}}$ in
characteristic three, which in turn leads to
$\tC(-1)^*=\SClif[0]{\bbF^2,\bb{2}}$ in characteristic three.

\nm{2} The previous remark does not completely solve the isomorphism
problem for the conclusions to Theorem \ref{thm-2gen-jordan-field},
since that theorem actually provides a classification up to
isomorphism of axial algebras of Jordan type $\eta$ equipped with two
marked generators. As mentioned under \ref{eg-jordan}, isomorphic
$2$-generated algebras $A$ and $A'$ can nevertheless give rise to
nonisomorphic marked algebras $(A,a,b)$ and $(A',a',b')$. Section 4 of
\cite{RS} discusses categories of marked algebras in detail.
\trats

Now consider a primitive axial algebra $A$ of Jordan type $\eta$ with
generating axis set $\axisset$ of arbitrary size (not necessarily
two).  Recall that $\bar{\axisset}$ is the smallest set of axes with
the properties: \startnm{i}
\item $\axisset \subseteq \bar{\axisset}$.
\item If $p \in \bar{\axisset}$ and $t$ is the Miyamoto involution
  associated with $p$, then $\bar{\axisset}^t \subseteq
  \bar{\axisset}$.
\tratsnm
Corollary \ref{cor-2gen} states that the algebra $A$ is the
$\bbF$-space spanned by the axes of $\bar{\axisset}$.

\dspace{.5}

\noindent
{\sc Proof of Corollary \ref{cor-2gen}.}

The algebra $A$ is spanned as $\bbF$-space by the multiplicative
submagma generated by $\axisset$.  By Theorem
\ref{thm-2gen-jordan-field} and Proposition \ref{prop-image-jordan}
every product of two members of $\axisset$, and indeed any two members
of $\bar{\axisset}$, is in the $\bbF$-span of
$\bar{\axisset}$. Therefore the span of $\bar{\axisset}$ is closed
under multiplication and contains the generators $\axisset$; it is
equal to $A$. \done

\section{Automorphisms of axial algebras of Jordan type}
\label{sec-auto}

In this section we focus on the automorphism groups of axial algebras
of Jordan type $\eta$ that are generated by Miyamoto
involutions. Especially we examine dihedral subgroups generated by two
Miyamoto involutions.

The following observation will be used frequently without mention.

\start{\lem}{lem-auto-action}
If $t$ is an automorphism of $A$ and $m$ is an axis, then $m^t$ is an
axis with $\tau(m)^t=\tau(m^t)$. \done
\trats

We first consider axes $a$ for which $\tau(a)=1$.

\start{\lem}{lem-tau-one-bis}
Let $A$ be an axial algebra of Jordan type $\eta$ over the field
$\bbF$ that is generated by the set $\axisset$ of axes.  Write
$\axisset$ as the disjoint union of $\axisset^1$ and $\axisset^\eta$,
where $a \in\axisset^1 $ if and only if $\tau(a)=1$ and $a \in
\axisset^\eta$ if and only if $\tau(a)$ has order $2$.  Then $A =
\left(\bigoplus_{a \in \axisset^1} \bbF a\right) \oplus A^\eta$ where
$A^\eta$ is the axial algebra of Jordan type $\eta$ generated by
$\axisset^\eta$.
\trats

\pf By Proposition \ref{prop-automorphism}, for an algebra $A$ of
Jordan type we have $\tau(a)=1$ precisely when $A=\bbF a \oplus
A_0(a)$. In this case by Proposition \ref{prop-0-1} the subalgebra
$A_0(a)$ contains all the axes of $\axisset$ except $a$. Especially
$\bigcap_{a \in \axisset^1} A_0(a)$ is the subalgebra $A^\eta$
generated by $\axisset^\eta$ and then $A = \left(\bigoplus_{a \in
  \axisset^1} \bbF a\right) \oplus A^\eta$, as claimed.  \done

\dspace{.5}

\noindent
Thus we are able to reduce to the case where all $\tau(m)$ have order
$2$ and are Miyamoto involutions.

For $\eta=\frac{1}{2}$, the examples coming from Clifford algebras
(and described in \ref{eg-cliff}) have
$\tau(e_0^+)\tau(e_1^+)=t_{v_0}t_{v_1}$, the product of two orthogonal
reflections on the space $\bbF^2$.  All orders are possible (although
restrictions on the field $\bbF$ and the form $b$ would lead to order
restrictions). Especially, for infinite $\bbF$ it is possible for the
two Miyamoto involutions to generate an infinite dihedral group so
that finite $\axisset$ (of size two) generates an a algebra of finite
dimension but with $\bar{\axisset}$ infinite.  We will have little
more to say regarding the case $\eta=\nfrac{1}{2}$ in this section.

For $\eta \neq\frac{1}{2}$, the situation is much different.  Theorem
\ref{thm-2gen-jordan-field} tells us that the possibilities for the
dihedral group generated by two Miyamoto involutions are extremely
limited.

The remainder of this section is focused on proving a more precise
version of Theorem \ref{thm-miyamoto-neq2}:

\start{\thm}{thm-miyamoto-neq2-bis}
Let $A$ be an axial algebra of Jordan type $\eta \neq \frac{1}{2}$
over a field of characteristic not two that is generated by the set
$\axisset$ of axes.

Write $\axisset$ as the disjoint union of $\axisset^1$ and
$\axisset^\eta$, where $a \in\axisset^1 $ if and only if $\tau(a)=1$
and $a \in \axisset^\eta$ if and only if $\tau(a)$ has order $2$.
\startnm{a}
\item
$A = \left(\bigoplus_{a \in \axisset^1} \bbF a\right) \oplus A^\eta$
  where $A^\eta$ is the axial algebra of Jordan type $\eta$ generated
  by $\axisset^\eta$.
\item
The map $a \mapsto \tau(a)$ is a bijection of $\bar{\axisset}^\eta$
with the corresponding set $D$ of Miyamoto involutions, and $D$ is a
normal set of $3$-transpositions in the subgroup $\langle D \rangle$
of the automorphism group of $A$ and $A^\eta$.
\tratsnm
\trats

A normal set $D$ of elements of order $2$ in the group $G$ is said to
consist of \alrt{$3$-transpositions} provided, for each pair $d,e \in
D$ the order of the product $de$ is $1$, $2$, or $3$.  Equivalently,
we must (respectively) have one of $d=e$, $de=ed\neq 1$, or $\langle
d,e \rangle \iso \Sym{3}$.

The elements of $D$ are called \alrt{transpositions} as the motivating
examples are the transpositions of any symmetric group. A group
generated by a normal set of $3$-transpositions is called a
\alrt{$3$-transposition group}.  These have been studied extensively
since their introduction by Fischer \cite{fischer}.  Fischer's work
and its successors, especially \cite{CuHa}, effectively classify all
$3$-transposition groups.

\start{\prop}{prop-2gen-miyamoto}
Let $A$ be an axial algebra of Jordan type $\eta$.  Let $a$ and $b$ be
two axes of $A$ with $\tau(a)$ and $\tau(b)$ the corresponding
Miyamoto involutions, and let $N$ be the subalgebra generated by $a$
and $b$.
\startnm{a}
\item If $N$ is of type $\oA$, then $a=b$, $\tau(a)=\tau(b)$, and
  $\tau(a)\tau(b)=1$.
\item If $N$ is of type $\tB$, then $\tau(a)\tau(b)=\tau(b)\tau(a)$
  and $(\tau(a)\tau(b))^2=1$.
\item If $N$ is of type $\tC(\eta)$ or type $\tC(-1)^*$, then
  $\tau(a)^{\tau(b)}=\tau(b)^{\tau(a)}$ and $(\tau(a)\tau(b))^3=1$.
\tratsnm
\trats

\pf
We consider the cases in turn.

\nm{a} Here $N$ contains a single axis $a=b$, so $\tau(a)=\tau(b)$.
 
\nm{b} By \ref{eg-tB} we have $b^{\tau(a)}=b$, so $\tau(b)^{\tau(a)}
=\tau(b)$ and $(\tau(a)\tau(b))^2=\tau(b)^{\tau(a)}\tau(b)=1$.

\nm{c}
By \ref{eg-tC} we have
\(
b^{\tau(a)}=c=a^{\tau(b)}\,.
\)
Therefore 
\[
\tau(b)^{\tau(a)}=\tau(c)=\tau(b)^{\tau(a)}
\]
and
\[
(\tau(a)\tau(b))^3=\tau(b)^{\tau(a)}\tau(a)^{\tau(b)}
=\tau(c)\tau(c)=1\,. \done
\]

\dspace{.5}

\start{\prop}{prop-tau-equal}
Let $A$ be a primitive axial algebra of Jordan type $\eta$ over a
field of characteristic not two that is generated by the set
$\axisset$ of axes.  Assume that every subalgebra generated by two
elements of $\axisset$ has type one of $\oA$, $\tB$, $\tC(\eta)$, or
$\tC(-1)^*$.  (By Theorem \ref{thm-2gen-jordan-field}, this is the
case for $\eta \neq \onehalf$.)

If $a,b \in \axisset$ with $\tau(a)=\tau(b)\neq 1$ then $a=b$.
\trats

\pf Let $t=\tau(a)=\tau(b)$ and choose $c\in \axisset$ with $c^t\neq
c$, possible as $t \neq 1$.  Let $N_{a,c}$ be the subalgebra generated
by $a$ and $c$, and let $N_{b,c}$ be the subalgebra generated by $b$
and $c$.  By hypothesis the types of $N_{a,c}$ and $N_{b,c}$ are
$\tC(\eta)$ or possibly $\tC(-1)^*$ (when $\eta=-1$). In particular,
$|\bar{\axisset} \cap N_{a,c}|= |\bar{\axisset} \cap N_{b,c}|=3$ with
$\langle \tau(a),\tau(c) \rangle = \langle \tau(b),\tau(c)\rangle$
acting as the symmetric group of degree $3$ on each.  But then
\[
a^{\tau(c)}=c^{\tau(a)}=c^t=c^{\tau(b)}=b^{\tau(c)}\,,
\]
hence $a=b$.\done

\dspace{.5}

By \ref{prop-clif-prop} when $\eta=\nfrac{1}{2}$ distinct axes can
have the same Miyamoto involution.

By Lemma \ref{lem-tau-one-bis} it is possible to have
$\tau(a)=\tau(b)= 1$ for distinct $a$ and $b$ in $\axisset$, but then
the corresponding $1$-dimensional subalgebras can be ``subtracted
out.''

\dspace{.5}

\noindent
{\sc Proof of Theorem \ref{thm-miyamoto-neq2-bis}.}

Part \nm{a} follows directly from Lemma \ref{lem-tau-one-bis}.

As $\eta\neq\onehalf$ by assumption, Proposition \ref{prop-tau-equal}
says that the map $a \mapsto \tau(a)$ is a bijection of
$\bar{\axisset}^\eta$ and $D=\set{\tau(a)}{a \in
  \bar{\axisset}^\eta}$.  By definition
$(\bar{\axisset}^\eta)^t=\bar{\axisset}^\eta$ for each $t\in D$, so
$D^t=D$ is a normal subset in the subgroup $\langle D \rangle$ of the
automorphism groups of $A$ and $A^\eta$.  Theorem
\ref{thm-2gen-jordan-field} tells us that any two elements $a$ and $b$
of $\bar{\axisset}^\eta$ generate a subalgebra of type $\oA$, $\tB$,
$\tC(\eta)$, or $\tC(-1)^*$. Then by Proposition
\ref{prop-2gen-miyamoto} the corresponding product $\tau(a)\tau(b)$ of
two elements from $D$ has order $1$, $2$, or $3$.  This gives \nm{b}.
\done

\dspace{1}

Theorem \ref{thm-miyamoto-neq2-bis} includes Theorem
\ref{thm-miyamoto-neq2}.  Corollary \ref{cor-neq2-loc-fin} then states
that an axial algebra $A$ of Jordan type $\eta \neq \frac{1}{2}$
generated by a finite number of axes are finite dimensional.

\dspace{.5}

\noindent
{\sc Proof of
Corollary \ref{cor-neq2-loc-fin}.}

Theorem \ref{thm-miyamoto-neq2-bis} allows us to assume that $\tau(a)$
has order $2$ for every $a \in \axisset$ and that there is a bijection
between the Miyamoto involutions of the normal $3$-transpositions set
$D$ and the elements of $\bar{\axisset}$.

By results from \cite{CuHa}, every finitely generated
$3$-transposition group is finite. In particular, the group generated
by the Miyamoto involutions for finite $\axisset$ is a finite group.
As that subgroup contains the Miyamoto involutions for
$\bar{\axisset}$, that set too is finite.  The algebra is the
$\bbF$-span of $\bar{\axisset}$ by Corollary \ref{cor-2gen}, so the
algebra $A$ is finite dimensional.  \done

\dspace{.5}

In fact, $|\bar{\axisset}|$ and so the dimension of $A$ can be bounded
by a function of $|\axisset|$; see \cite{HaSh}.

\section{Matsuo algebras and Fischer spaces}
\label{sec-matsuo}

Versions of certain results from this section appeared originally in
the unpublished work of Matsuo and Matsuo, \cite{MaMa} and
\cite{Matsuo-ver1}. Related results also appear in \cite{Re13}.

A \alrt{partial triple system} or \alrt{partial linear space of order
  two} $\Pi=(\calP,\calL)$ is a set $\calP$, called \alrt{points}, and
a set $\calL$ of subsets of $\calP$, called \alrt{lines}, such that:
\startnm{i}
\item
every line contains exactly three points;
\item
two distinct points belong to at most one line.
\tratsnm
If every pair of distinct points belongs to a unique line, then we
have a \alrt{linear space of order two} or \alrt{Steiner triple
  system}.

For distinct points $p,q$, we write $p \sim q$ if $p$ and $q$ are
collinear and $p \perp q$ if they are not. The set $p^\perp$ is then
the set of all points not collinear with $p$.  Let $\approx$ be the
equivalence relation on $\calP$ generated by $\sim$. The
$\approx$-equivalence classes are the \alrt{connected components} of
$\Pi$.

A \alrt{subspace} $(\calP',\calL')$ of $(\calP,\calL)$ is a subset
$\calP'$ of $\calP$ with the property that whenever there is a line of
$\calL$ intersecting $\calP'$ in at least two points, then the line is
a subset of $\calP'$ and so a member of the line set $\calL' \subseteq
\calL$.  For instance, every line is itself a subspace. A \alrt{plane}
is the subspace generated by two distinct, intersecting lines---the
intersection of all subspaces containing the two lines. Each connected
component of $\Pi$ is a subspace, and $\Pi$ is then the disjoint union
of its connected components.

Of particular interest here are the \alrt{Fischer spaces}.  If $D$ is
a normal set of $3$-transpositions in the group $G$, then the
associated Fischer space $\Pi$ is (up to isomorphism) the partial
triple system having point set $D$ and line set consisting of the
triples of points (transpositions) from the various subgroups
$\Sym{3}$ generated by two transpositions.

For each $p \in \calP$ the associated transposition $t(p)$ of $D$ is
an automorphism of $\Pi$ that fixes $p$ and each $q$ of $p^\perp$
while switching $r$ and $s$ whenever $\{p,r,s\}$ is a line. In
particular the subspace fixed by $t(p)$ is the disjoint union of
$\{p\}$ and the subspace $p^\perp$.

Fischer spaces have a well known characterization due to Buekenhout:

\start{\prop}{prop-buek}
A partial triple system is a Fischer space if and only if the only
isomorphism types allowed for its planes are the dual affine plane of
order two and the affine plane of order three.  \done
\trats

Let $\Pi=(\calP,\calL)$ be a partial triple system.  Choose a constant
$\delta$ in the field $\bbF$.  The associated \alrt{Matsuo algebra}
$\Matsuo{\Pi,\delta,\bbF}$ over the field $\bbF$ is the $\bbF$-space
$\bigoplus_{p \in \calP} \bbF a_p$ with multiplication provided by:
\startnm{i}
\item
For $p \in \calP$ we have $a_p^2=a_p$.
\item 
For distinct $p,q \in \calP$ with $p$ and $q$ not collinear,
$a_pa_q=0$.
\item 
For distinct $p,r \in \calP$ with $\{p,r,s\}$ a line,
$a_pa_r=\delta(a_p + a_r - a_s)$.
\tratsnm

For each line $\{p,r,s\}$ set
\[
a_{prs}= \eta a_p -a_r-a_s\quad \text{and}\quad g_{prs}= a_r-a_s\,.
\]

\start{\thm}{thm-conn-comp}
The Matsuo algebra $\Matsuo{\Pi,\nfrac{\eta}{2},\bbF}$ is a primitive
axial algebra for the basis $\axisset=\set{a_p}{p\in \calP}$ of axes
with eigenvalue set $\{1,0,\eta\}$.  It is the direct sum of its
ideals $\Matsuo{\Pi^{(i)},\nfrac{\eta}{2},\bbF}$, where the
$\Pi^{(i)}$, $i \in I$, are the connected components of $\Pi$.
\trats

\pf For each $p \in \calP$ and each line $\{p,r,s\}$ the pair
$\{a_r,a_s\}$ from the canonical basis $\axisset$ can be replaced by
$\{a_{prs},g_{prs}\}$.  (While $a_{prs}=a_{psr}$, we choose only one
of $g_{prs}$ and $g_{psr}=-g_{prs}$.)  Then by \ref{eg-tC} we have
\begin{align*}
M_1(a_p) &= \bbF a_p\\
M_0(a_p) 
&= \bigoplus_{q \in p^\perp} \bbF a_q\
\oplus 
\bigoplus_{\{p,r,s\} \in \calL} \bbF\, a_{prs}
\\
M_\eta(a_p) &= 
\bigoplus_{\{p,r,s\} \in \calL} \bbF g_{prs}\,.
\end{align*}
In particular $M$ has axial basis $\axisset=\set{a_p}{p \in \calP}$
with eigenvalue set $\{1,0,\eta\}$, as claimed.

If $p$ and $q$ are points in different connected components of $\Pi$,
then $a_pa_q=0$. Therefore each $\axisset^{(i)}
=\set{a_p}{p\in\Pi^{(i)}}$ is the basis of a subalgebra $M^{(i)}$
isomorphic to $\Matsuo{\Pi^{(i)},\nfrac{\eta}{2},\bbF}$ and such that
$M^{(i)}M^{(j)}=0$ for $i \neq j$.  Especially each $M^{(i)}$ is an
ideal, and the algebra is the direct sum of these ideals.  \done

\dspace{1}

From Theorem \ref{thm-miyamoto-neq2-bis} we have immediately

\start{\thm}{thm-matsuo-exist}
For $\eta \neq \onehalf$, every primitive axial algebra of Jordan type
$\eta$ is isomorphic to $(\bigoplus_{i \in I} \bbF)\oplus M$, for some
index set $I$, where $M$ is a quotient of the Matsuo algebra
$\Matsuo{\Pi,\nfrac{\eta}{2},\bbF}$ associated with the Fischer space
$\Pi$ of all Miyamoto involutions. \done
\trats

The next theorem gives a construction that provides a converse but
also works in the case $\eta = \onehalf$.

\start{\thm}{thm-matsuo}
Let $\bbF$ be a field of characteristic not two and $\eta \in \bbF$
with $\eta \neq 0,1$.  If $\Pi=(\calP,\calL)$ is a Fischer space, then
the associated Matsuo algebra $M=\Matsuo{\Pi,\nfrac{\eta}{2},\bbF}$ is
a primitive axial algebra of Jordan type $\eta$ generated by its basis
of axes $\axisset=\set{a_p}{p \in \calP}$.  Each transposition $t(p)$
acts as the Miyamoto involution $\tau(a_p)$.
\trats

\pf By Theorem \ref{thm-conn-comp} the algebra $M$ is primitive and
axial with the basis $\axisset$ of axes for the eigenvalues
$\{1,0,\eta\}$.

By \ref{eg-tC} (and as observed in the proof of Theorem
\ref{thm-conn-comp}) $M_0(a_p)$ is spanned by the elements $a_q$, for
$q \in p^\perp$, and $a_{prs}= \eta a_p -a_r-a_s$, for $\{p,r,s\} \in
\calL$, while $M_\eta(a_p)$ is spanned by the elements $g_{prs}=
a_r-a_s$, for $\{p,r,s\} \in \calL$.

The transposition $t(p)$ acts on the basis elements of $\axisset$ via
$a_q^{t(p)}=a_{q^{t(p)}}$. Thus, in view of the previous paragraph,
$t(p)$ induces on $M$ the Miyamoto involution $\tau(a_p)$ for the
$\zt$-grading $M_+(a_p)=M_1(a_p)\cup M_0(a_p)$ and
$M_-(a_p)=M_\eta(a_p)$.

By Lemma \ref{lem-subalg-bis}, it remains to check that each
$M_0(a_p)$ is a subalgebra---that it is closed under multiplication by
its spanning elements.  Thus there are three cases to consider.
\startnm{1}
\item {\em $a_qa_v$ with $q,v \in p^\perp$.}

Certainly $a_q^2=a_q$.  If $q \in v^\perp$ then $a_qa_v=0 \in
M_0(a_p)$.  If $\{q,v,w\} \in \calL$, then
$a_qa_v=\halfeta(a_q+a_v-a_w)$.  Since $p^\perp$ is a subspace
containing $q$ and $v$, it also contains $w$. Thus $a_w$ and $a_qa_v$
are both in $M_0(a_p)$.

\item {\em $a_qa_{prs}$ with $q\in p^\perp$ and $\{p,r,s\}\in\calL$.}

If $r$ or $s$ is in the subspace $q^\perp$, then the entire line
$\{p,r,s\}$ is in $q^\perp$ and $a_qa_{prs}=0\in M_0(a_p)$.

Now assume that $r,s \notin q^\perp$. Then $p,q,r,s$ all belong to the
plane generated by $\{p,r,s\}$ and the line through $q$ and $r$. As $q
\in p^\perp$, this plane must be dual affine of order two, and we may
take its line set to be
\[
\{p,r,s\}\,,\ \{p,t,u\}\,,\ \{q,r,t\}\,,\ \{q,s,u\}\,.
\]
Thus 
\begin{align*}
a_qa_{prs} &= a_q(\eta a_p-a_r-a_s)\\
&=\eta a_pa_q-a_qa_r-a_qa_s\\
&=-\halfeta(a_q+a_r-a_t)-\halfeta(a_q+a_s-a_u)\\
&=-\eta a_q+\halfeta(-a_r-a_s)-\halfeta(-a_t-a_u)\\
&=-\eta a_q+\halfeta a_{prs}-\halfeta a_{ptu} \in M_0(a_p)\,.
\end{align*}

\item {\em $a_{prs}a_{ptu}$ with $\{p,r,s\},\{p,t,u\}\in\calL$.}

We have $a_{prs}^2 \in M_0(a_p)$ by \ref{eg-tC}.  For two distinct
intersecting lines $\{p,r,s\}$ and $\{p,t,u\}$, the calculation
ultimately depends upon the type of the plane $\Delta$ they
generate. In either case
\begin{align*}
a_{prs}a_{ptu}
&=(\eta a_p -a_r -a_s)(\eta a_p -a_t -a_u)\\
&=(-a_r -a_s) (\eta a_p -a_t -a_u)\\
&=-\eta a_ra_p -\eta a_sa_p 
+a_ra_t+a_ra_u+a_sa_t+a_sa_u\\
&=-\nfrac{\eta^2}{2}(a_p+a_r-a_s)  -\nfrac{\eta^2}{2}(a_p+a_s-a_r)\\
&\qquad +a_ra_t+a_ra_u+a_sa_t+a_sa_u\\
&= -\eta^2 a_p + a_ra_t+a_ra_u+a_sa_t+a_sa_u\,.
\end{align*}
If $\Delta$ is dual affine of order two, then we may take its lines to
be those of \nm{2}, so that $a_ra_u=0=a_sa_t$ and
\begin{align*}
a_{prs}a_{ptu}
&= -\eta^2 a_p + a_ra_t+a_ra_u+a_sa_t+a_sa_u\\
&= -\eta^2 a_p + a_ra_t+a_sa_u\\
&= -\eta^2 a_p +\halfeta(a_r+a_t-a_q)  +\halfeta(a_s+a_u-a_q)\\
&= -\eta a_q -\halfeta(2\eta\, a_p -a_r-a_t-a_s-a_u)\\
&= -\eta a_q -\halfeta(a_{prs} + a_{ptu})\in M_0(a_p)\,.
\end{align*}

On the other hand, if $\Delta$ is affine of order three, then we may
take the two additional lines of $\Delta$ on $p$ to be $\{p,w,x\}$ and
$\{p,y,z\}$. We then have
\begin{align*}
a_{prs}a_{ptu}
&= -\eta^2 a_p + a_ra_t+a_ra_u+a_sa_t+a_sa_u\\
&= -\eta^2 a_p +\eta(a_r+a_s+a_t+a_u)  -\halfeta(a_w+a_x+a_y+a_z)\\
&=  -\halfeta\left( 2\eta\, a_p-2(a_r+a_s)-2(a_t+a_u) 
+(a_w+a_x)+(a_y+a_z)\right)\\
&= -\halfeta(2a_{prs} + 2a_{ptu} -a_{pwx}-a_{pyz})
\in M_0(a_p)\,. \done
\end{align*}
\tratsnm

Theorem \ref{thm-matsuo} gives all parts of Theorem
\ref{thm-matsuo-bis} of the introduction except for the existence of
an associative form, which is handled in Corollary
\ref{cor-matsuo-frob} below.

\dspace{.5}

We have a pleasant consequence of the work in this section:

\start{\thm}{thm-felix}
The Matsuo algebra $\Matsuo{\Pi,\delta,\bbF}$ is an axial algebra of
Jordan type if and only if $\Pi$ is a Fischer space.
\trats

\pf
The only possible Jordan type is $\eta=2\delta$.

Theorem \ref{thm-matsuo} gives the converse part of this theorem
immediately.  For $\eta \neq \onehalf$ the rest follows from Theorem
\ref{thm-matsuo-exist}, but this difficult result is not necessary in
proving the direct part for arbitrary $\eta$.

Suppose that the Matsuo algebra $\Matsuo{\Pi,\delta,\bbF}$ is an axial
algebra of Jordan type $\eta=2\delta$ presented using the partial
triple system $\Pi=(\calP,\calL)$. Then for each $x \in \calP$, the
Miyamoto involution $\tau(a_x)$ permutes the generating set
$\axisset_\calP=\set{a_p}{p \in \calP}$, taking $\tB$ subalgebras to
$\tB$ subalgebras and $\tC(\eta)$ algebras to $\tC(\eta)$ algebras.
Therefore the induced permutation $t(x)$ of $\calP$ given by
$a_p^{\tau(a_x)}=a_{p^{t(x)}}$ is an automorphism of the partial
triple system $\Pi$. Indeed it is the unique \alrt{central}
automorphism of $\Pi$ with \alrt{center} $x$---that is, it fixes $x$
and all points not collinear with $x$ and, for each line $\{x,y,z\}$
on $x$, it switches $y$ and $z$.

It is well-known, and easy to check, that a triple system admits all
possible central automorphisms if and only if the collection of
central automorphisms is a normal set of $3$-transpositions in
$\Aut{\Pi}$ with $\Pi$ as the corresponding Fischer space. Indeed, for
any automorphism $g$ we always have $t(x)^g=t(x^g)$. In particular if
distinct $x$ and $y$ are not collinear, then
\[
t(y)^{t(x)}=t(x)t(y)t(x)=t(y)\,,
\]
and $(t(x)t(y))^2=1$,
while if $\{x,y,z\}$ is a line 
\[
t(x)t(y)t(x)=t(y)^{t(x)}=t(z)=t(x)^{t(y)}=t(y)t(x)t(y)
\]
and $(t(x)t(y))^3=t(z)^2=1$. Therefore $\Pi$ is a Fischer space. \done

\section{Frobenius axial algebras of Jordan type}
\label{sec-frob}

The results from this section essentially appeared in
\cite{Matsuo-ver1}, the unpublished first version of the published
\cite{Matsuo}.

%
%
\newpage
\start{\lem}{lem-line-frob}
\startnm{a}
\item
The algebra $\tB$ is a Frobenius algebra.  The bilinear form $\lla
\cdot,\cdot \rra$ on this algebra is associative if and only if $ \lla
b_0,b_1 \rra = 0 = \lla b_1,b_0 \rra $.
\item
The algebra $\tC(\eta)$ is a Frobenius algebra.  An associative
bilinear form on this algebra is a scalar multiple of the form given
by
\[
\lla c_i,c_i \rra = 1\,,\
\lla c_i,c_j \rra = \nfrac{\eta}{2}
\ \text{for}\ i \neq j
\,.
\]
\tratsnm
\trats

\pf By Proposition \ref{prop-frob}(a) an associative form on these
axial algebras is symmetric.

\nm{a} The given forms are clearly associative. Now consider an
arbitrary associative form. Then for $i \neq j$,
\[
\lla b_i,b_j\rra=\lla b_ib_i,b_j\rra
=\lla b_i,b_ib_j\rra=\lla b_i,0\rra=0\,.
\]

\nm{b} For any bilinear form $\lla \cdot,\cdot\rra$ on $\tC(\eta)$ we
have, for $\{i,j,k\}=\{0,1,2\}$: \startnm{i}
\item 
\qquad\qquad\qquad
$\lla c_i,c_ic_i \rra= \lla c_i,c_i \rra =\lla c_ic_i,c_i \rra$ ;
\item 
\begin{align*}
\lla c_i,c_ic_j \rra
&= \halfeta\lla c_i,c_i+c_j-c_k \rra\\
&= \halfeta\left(\lla c_i,c_i\rra+\lla c_i,c_j\rra
-\lla c_i,c_k \rra\right)\,;
\end{align*}
\item 
\begin{align*}
\lla c_i,c_jc_k \rra - \lla c_ic_j,c_k\rra
&= \halfeta\lla c_i,c_j+c_k-c_i \rra-\halfeta\lla c_i+c_j-c_k,c_k \rra
\\
&= \halfeta\left(\lla c_i,c_j\rra-\lla c_j,c_k\rra-\lla c_i,c_i \rra
+\lla c_k,c_k \rra\right)\,.
\end{align*}
\tratsnm Therefore the form is associative if and only if the
righthand side of \nm{ii} is always equal to $\lla c_i,c_j\rra$ and
the righthand side of \nm{iii} is always equal to $0$. In particular,
the given values do lead to an associative form (as promised earlier
under \ref{eg-tC}).

Now assume that the form is associative.  By \ref{eg-tC},
$\bbF(c_j-c_k)$ is the $\eta$-eigenspace for $\ad_{c_i}$; so by
Proposition \ref{prop-frob} always $\lla c_i,c_j-c_k\rra=0$.
Therefore
\[
\lla c_0,c_1 \rra=
\lla c_0,c_2 \rra=
\lla c_1,c_2 \rra= \kappa\,,
\]
for some constant $\kappa$. Then as the righthand side of \nm{iii} is
always $0$,
\[
\lla c_0,c_0 \rra=
\lla c_1,c_1 \rra=
\lla c_2,c_2 \rra=k
\]
is constant as well.

Finally as the form is associative, \nm{ii} becomes
\[
\kappa=\lla c_i,c_j\rra=\lla c_i,c_ic_j\rra=
\halfeta\left(\lla c_i,c_i\rra+\lla c_i,c_j\rra-\lla c_i,c_k \rra\right)=
\halfeta k\,.
\] 
The constant $k$ determines the associative form up to a scalar
multiple, as claimed, the given values corresponding to $k=1$.  \done

\start{\thm}{thm-matsuo-frob}
Let $\Pi =(\calP,\calL)$ be a partial triple system.  The Matsuo
algebra $\Matsuo{\Pi,\nfrac{\eta}{2},\bbF}$ admits a nonzero
associative form if and only if
\startnm{i}
\item for each $x \in \calP$, the subset $x^\perp$ is a subspace of
  $\Pi$, and
\item if $\{x,y,z\}$ and $\ell=\{x,v,w\}$ are lines of $\calL$, then
  $\ell \cap y^\perp = \emptyset$ if and only if $\ell\cap
  z^\perp=\emptyset$.
\tratsnm 
When $\Pi$ is connected, such a form is a scalar multiple of the form
given by, for distinct $p,q \in \calP$,
\[
\lla a_p,a_p \rra =1\,;\
\lla a_p,a_q \rra =0\
\text{if}\
q \in p^\perp\,;\
\lla a_p,a_q \rra =\halfeta\
\text{if}\
q \notin p^\perp\,.
\]
\trats

\pf By the previous lemma, under any associative form the ideals
corresponding to the distinct connected components of $\Pi$ are
perpendicular. Furthermore, the form when restricted to a specific
component $\Pi_i$ can only be a scalar multiple of the given form.

It remains to prove that this does give an associative form for the
(connected) space $\Pi\,(=\Pi_i)$ if and only if the two conditions
\nm{i} and \nm{ii} are satisfied.

By linearity it suffices to prove
\[
\lla a_ra_p,a_t\rra = \lla a_r,a_pa_t \rra
\]
for all $r,p,t \in \calP$, where by Lemma \ref{lem-line-frob} we may
assume that $r,p,t$ are distinct and do not lie together in a line of
$\calL$.

If $r,t \in p^\perp$, then 
\[
\lla a_ra_p,a_t\rra =\lla 0,a_t\rra
=0= 
\lla a_r,0\rra=
\lla a_r,a_pa_t \rra\,,
\]
as desired. Therefore we may also assume that $r \sim p$.  Let
$\{p,r,s\}$ be the line on $r$ and $p$.

\begin{quote}
{\sc Claim.} {\em $\lla a_ra_p,a_t\rra = \lla a_r,a_pa_t \rra$ for all
  triples of points with $ r \sim p \perp t$ if and only if $x^\perp$
  is a subspace for all $x \in \calP$.}
\end{quote}

We have 
\begin{align*}
\lla a_ra_p,a_t\rra - \lla a_r,a_pa_t \rra
&=\halfeta\lla a_r + a_p - a_s,a_t \rra - \lla a_r, 0 \rra\\
&=\halfeta\left(\lla a_r,a_t\rra +\lla a_p,a_t \rra - \lla a_s,a_t \rra\right)\\
&=\halfeta\left(\lla a_r,a_t\rra - \lla a_s,a_t \rra\right)\,.
\end{align*}
For this to be $0$ we must have either $r \perp t \perp s$ or $r \sim
t \sim s$.  As $t \perp p$, this says that $t^\perp$ either contains
all of the line $\{p,r,s\}$ or it only contains $p$. This happens for
all lines $\{p,r,s\}$ with $p \in t^\perp$ precisely when $t^\perp$ is
a subspace.  This gives the claim.

\begin{quote}
{\sc Claim.}  {\em Assume that $x^\perp$ is a subspace for all $x \in
  \calP$.  Then $\lla a_ra_p,a_t\rra = \lla a_r,a_pa_t \rra$ for all
  triples of points with $ r \sim p \sim t$ if and only if, for all
  $\{x,y,z\}$ and $\ell=\{x,v,w\}$ lines of $\calL$, we have $y^\perp
  \cap \ell = \emptyset$ if and only if $z^\perp \cap
  \ell=\emptyset$.}
\end{quote}

Let $\{p,t,u\}$ be a line. As $\lla a_p,a_t \rra = \halfeta = \lla
a_r,a_p\rra$,
\begin{align*}
\lla a_ra_p,a_t\rra - \lla a_r,a_pa_t \rra
&=\halfeta\left(\lla a_r + a_p - a_s,a_t \rra 
- \lla a_r, a_p+a_t-a_u \rra\right)\\
&=\halfeta\left(\lla a_r,a_t\rra +\lla a_p,a_t \rra 
- \lla a_s,a_t \rra\right.\\
&\qquad\qquad \left.-\lla a_r,a_p\rra -\lla a_r,a_t \rra 
+ \lla a_r,a_u \rra\right)\\
&=\halfeta\left(- \lla a_s,a_t \rra+ \lla a_r,a_u \rra\right)\,.
\end{align*}
This is $0$ when $\lla a_r,a_u \rra=\lla a_s,a_t \rra $;
that is, when we have
\[
(*)\quad\
\text{either}\quad r \sim u\ \text{and}\  s \sim t\ 
\quad\text{or}\quad 
r \perp u\ \text{and}\ s \perp t\,.
\]
We also have $r \sim p \sim u$, so we may replace $t$ by $u$
in the above to find that for the form to be associative
on these two lines we must additionally have
\[
(**)\quad
\text{either}\quad r \sim t\ \text{and}\  s \sim u\ 
\quad\text{or}\quad 
r \perp t\ \text{and}\ s \perp u\,.
\]
Conversely, the validity of $(*)$ and $(**)$ is sufficient for the
form to be associative on the two lines.

As all $q^\perp$ for $q \in \{r,s,t,u\}$ are subspaces, each of
$\{r,s\}$ must be collinear with at least one of $\{t,u\}$ and vice
versa. This is equivalent to $(*)$ and $(**)$ except for the
possibility that one of $r$ and $s$ is collinear with both of $t$ and
$u$ while the other is collinear with only one.  To avoid this, we
want $r$ to be collinear with both $t$ and $u$ if and only if $s$ is
as well.  This is equivalent to requiring that $r^\perp \cap
\{p,t,u\}$ is empty if and only if $s^\perp \cap \{p,t,u\}$ is empty.

This completes our proof of the second claim 
and so of the theorem. \done

\start{\cor}{cor-steiner}
If $\Pi$ is a Steiner triple system, then
$\Matsuo{\Pi,\nfrac{\eta}{2},\bbF}$ admits an associative form, which
is uniquely determined up to a scalar multiple as the form given by,
for distinct $p,q \in \calP$,
\[
\lla a_p,a_p \rra =1\,;\
\lla a_p,a_q \rra =0\
\text{if}\
q \in p^\perp\,;\
\lla a_p,a_q \rra =\halfeta\
\text{if}\
q \notin p^\perp\,.
\]
\trats

\pf In Steiner triple systems, each $x^\perp$ is empty. \done

\start{\cor}{cor-matsuo-frob}
If $\Pi$ is a Fischer space, then $\Matsuo{\Pi,\nfrac{\eta}{2},\bbF}$
is a Frobenius axial algebra of Jordan type $\eta$.  When $\Pi$ is
connected, an associative form is uniquely determined up to a scalar
multiple as the form given by, for distinct $p,q \in \calP$,
\[
\lla a_p,a_p \rra =1\,;\
\lla a_p,a_q \rra =0\
\text{if}\
q \in p^\perp\,;\
\lla a_p,a_q \rra =\halfeta\
\text{if}\
q \notin p^\perp\,.
\]
\trats

\pf The fixed point subspace of the transposition $t(x)$ is $\{x\}
\cup x^\perp$, the disjoint union of $\{x\}$ and the subspace
$x^\perp$, giving \nm{i} of Theorem
\ref{thm-matsuo-frob}. Furthermore, if $\{x,y,z\}$ and $\ell$ are two
lines on the point $x$, then $t(x)$ switches $y$ and $z$ and leaves
$\ell$ fixed globally, so that $(y^\perp\cap\ell)^{t(x)}=
z^\perp\cap\ell$, giving \nm{ii} of the theorem.  \done

\dspace{.5}

This corollary has at least three uses. Especially, it completes our
proof of Theorem \ref{thm-miyamoto-neq2}.

By Proposition \ref{prop-frob} and Theorem \ref{thm-conn-comp} every
ideal of the algebra is a sum of ideals corresponding to connected
components of $\Pi$ and ideals contained in the radical, which is the
maximal ideal containing no axes.

Finally, in the traditional applications the algebra is defined over
$\bbR$ and comes equipped with an associative, positive definite
form. The form $\lla \cdot,\cdot \rra$ of Corollary
\ref{cor-matsuo-frob} has Gram matrix $I + \nfrac{\eta}{2} D$, where
$D$ is the adjacency matrix of the $\Pi$-collinearity graph on
$\calP$.  In particular the form is positive (semi)definite when $D$
has minimal eigenvalue greater than (or equal to) $-2\eta^{-1}$.  For
Griess and Majorana algebras, the cases of interest are
$\eta=\nfrac{1}{4}$ and $\eta=\nfrac{1}{32}$.  These give the minimal
eigenvalues $-8$ (studied by Matsuo \cite{Matsuo}) and $-64$
(considered by Hall and Shpectorov \cite{HaSh}).

\end{document}